\documentclass[11pt]{article}
\usepackage[left=1in,right=1in,top=1in,bottom=1in]{geometry}
\usepackage{times}
\usepackage{expl3}
\usepackage{cite}
\usepackage[table]{xcolor}
\usepackage{multirow}
\usepackage{stackengine} 
\usepackage{hhline}
\usepackage{lipsum}
\usepackage{titlesec}
\usepackage{wrapfig}
\usepackage{epsfig}
\usepackage{bbm}
\usepackage{graphicx}
\usepackage{amsmath}
\usepackage[title]{appendix}
\usepackage{amssymb}
\usepackage{epstopdf}
\usepackage{boldline}
\usepackage{calligra}
\usepackage{bm}
\usepackage{url}
\usepackage{blindtext}

\newcommand{\define}{\stackrel{\mbox{\tiny def}}{=}}

\newtheorem{definition}{Definition}
\newtheorem{theorem}{Theorem}
\newtheorem{corollary}{Corollary}
\newtheorem{lemma}{Lemma}

\usepackage{mathtools}
\usepackage{epstopdf}
\usepackage{balance}
\usepackage{thmtools}
\usepackage{thm-restate}
\usepackage{hyperref}
\usepackage{cleveref}

\usepackage[ruled,vlined]{algorithm2e}
\include{pythonlisting}
\newcommand{\ostar}{\mathbin{\mathpalette\make@circled\star}}

\makeatletter
\newcommand{\removelatexerror}{\let\@latex@error\@gobble}
\makeatother
\makeatletter
\newcommand*{\rom}[1]{\expandafter\@slowromancap\romannumeral #1@}
\makeatother

\ExplSyntaxOn
\newcommand\latinabbrev[1]{
  \peek_meaning:NTF . {% Same as \@ifnextchar
    #1\@}%
  { \peek_catcode:NTF a {% Check whether next char has same catcode as \'a, i.e., is a letter
      #1.\@ }%
    {#1.\@}}}
\ExplSyntaxOff

\titleclass{\subsubsubsection}{straight}[\subsubsection]

\begin{document}
\vspace{1cm}
\title{Random Tensor Inequalities and Tail bounds for Bivariate Random Tensor Means, Part I}\vspace{1.8cm}
\author{Shih~Yu~Chang 
% <-this % stops a space
\thanks{Shih Yu Chang is with the Department of Applied Data Science,
San Jose State University, San Jose, CA, U. S. A. (e-mail: {\tt
shihyu.chang@sjsu.edu}).
           }}

\maketitle

\begin{abstract}
In this work, we apply the concept about \emph{operator connection} to consider bivariate random tensor means. We first extend classical Markov's and Chebyshev’s inequalities from a random variable to a random tensor by establishing Markov's inequality for tensors and Chebyshev’s inequality for tensors. These inequalities are applied to establish tail bounds for bivariate random tensor means represented by operator perspectives based on various types of connection functions: tensor increasing functions, tensor decreasing functions, and tensor concavity functions. We also consider tail bounds relations for the summation and product of eigenvalues based on majorization ordering of eigenvalues of bivariate random tensor means. This is Part I of our work about random tensor inequalities and tail bounds for bivariate random tensor mean. In our Part II, we will consider bivariate random tensor mean with respect to non-invertible random tensors and their applications.  
\end{abstract}

\begin{keywords}
Markov's inequality, Chebyshev’s inequality, random tensors, bivariate tensor mean, L\"owner ordering, 
majorization ordering.
\end{keywords}

\section{Introduction}\label{sec:Introduction}

Random tensors have both theoretical and practical applications in various fields such as machine learning, physics, and computer science. In machine learning, random tensors are commonly used for weight initialization in neural networks. A neural network is a collection of interconnected nodes (neurons) that take input data, perform computations, and produce an output. The weights of these connections between neurons are typically initialized with random values before training. This helps to break the symmetry of the network and ensure that each neuron learns different features of the data~\cite{ouerfelli2022random}. In physics, random tensors are used in the study of quantum entanglement and the geometry of entangled states. In particular, random tensor networks have been used to simulate quantum systems and understand the properties of entangled states. These studies have applications in quantum computing, quantum field theory, and condensed matter physics~\cite{evnin2021melonic}. In computer science, random tensors are used in the design and analysis of algorithms. For example, randomized matrix algorithms use random tensors to efficiently compute matrix decompositions, which have applications in data analysis, signal processing, and machine learning. Randomized algorithms are also used in graph theory to solve problems such as graph partitioning and clustering~\cite{tropp2015introduction}. Overall, the theory and application of random tensors have a wide range of practical uses across different fields, from improving the performance of machine learning models to advancing our understanding of quantum physics and optimizing algorithms~\cite{chang2022TWF,chang2022randouble,chang2022tPII,chang2022tPI,
chang2022generaltail,chang2022TKF,chang2022tensorq}.

In our recent work, we studied bivariate random tensor means by the format of random double tensor integrals (DTI)~\cite{chang2022randouble}. The tail bound of the unitarily invariant norm for the random DTI is derived and this bound enable us to derive tail bounds of the unitarily invariant norm for various types of two tensors means, e.g., arithmetic mean, geometric mean, and harmonic mean. In this work, we apply the notion of \emph{operator connection} to explore bivariate random tensor means~\cite{kubo1980means}. Ando-Hiai type inequalities and their applications have attracted active research in the community of operator theory since 1994~\cite{ando1994log, wada2018does}, and they play a crucial role in recent evolution of bivariate/multivariate operator means. By treating a tensor as a operator, an \emph{operator perspective} of two tensors $\mathcal{A} \in \mathbb{C}^{I_1 \times \dots \times I_N\times I_1 \times \dots \times I_N}$ and $\mathcal{B} \in \mathbb{C}^{I_1 \times \dots \times I_N\times I_1 \times \dots \times I_N}$ is a two-argument tensor functions associated to a continuous function (A.K.A. connection function) $g$  on $(0, \infty)$, denoted by $\mathcal{A} \#_{g} \mathcal{B}$, which is defined by  as 
\begin{eqnarray}
\mathcal{A} \#_{g} \mathcal{B} \define \mathcal{B}^{1/2}\star_N g\left(\mathcal{B}^{-1/2}\star_N \mathcal{A} \star_N \mathcal{B}^{-1/2} \right)\star_N \mathcal{B}^{1/2}, 
\end{eqnarray}
where $\star_N$ is an Einstein product between two tensors defined by Eq.~\eqref{eq: Einstein product definition}~\cite{hiai2020ando,effros2014non}. When the function $g$ is a positive operator monotone function with $g(1)=1$, the operator perspective $\#_{g}$ operation is reduced as an operator mean operation~\cite{kubo1980means}.

In this work, we first extend classical Markov's and Chebyshev’s inequalities from a random variable to a random tensor by establishing Markov's inequality for tensors and Chebyshev’s inequality for tensors.  These inequalities are used to provide tail bounds for bivariate random tensor means represented by operator perspectives based on various types of connection functions: tensor increasing functions, tensor decreasing functions, and tensor concavity functions. We also derive tail bounds relations for the summation and product of eigenvalues based on majorization ordering of eigenvalues of bivariate random tensor means. Our work about the part of L\"owner ordering between bivariate random tensor means is based on recent work from~\cite{hiai2020ando}, however, several L\"owner ordering relations have been generalized to a larger range of exponent via recursion technique and Kantorovich type inequality. This is Part I of our work about random tensor inequalities and tail bounds for bivariate random tensor mean. In our Part II, we will consider bivariate random tensor mean with respect to non-invertible random tensors and their applications.  

The rest of this paper is organized as follows. In Section~\ref{sec:Tensor Basics}, we will review basic definitions about tensors. Tensor Markov's and tensor Chebyshev’s inequalities are present in Section~\ref{sec:Tensor Markov Chebyshev’s Inequalities}. In Section~\ref{sec:Tail Bounds for Bivariate Random Tensors Mean Based on Loewner Ordering}, we will establish several tail bounds for bivariate random tensor means with respect to various types of connection functions. In Section~\ref{sec:Tail Bounds for Bivariate Random Tensors Mean Based on Majorization Ordering},  we will explore tail bounds relations for the summation and product of eigenvalues based on majorization ordering of eigenvalues of bivariate random tensor means.

\noindent \textbf{Nomenclature:} The sets of complex and real numbers are denoted by $\mathbb{C}$ and $\mathbb{R}$, respectively. The set of natural numbers is represented by $\mathbb{N}$. A scalar is denoted by an either italicized or Greek alphabet such as $x$ or $\beta$; a vector is denoted by a lowercase bold-faced alphabet such as $\bm{x}$; a matrix is denoted by an uppercase bold-faced alphabet such as $\bm{X}$; a tensor is denoted by a calligraphic alphabet such as $\mathcal{x}$ or $\mathcal{X}$. 

\section{Tensor Basics}\label{sec:Tensor Basics}

In this section, we will review necessary basic facts about tensors. Given $\mathcal{X} \define (x_{i_1, \dots, i_M, j_1, \dots,j_N}) \in \mathbb{C}^{I_1 \times \dots \times I_M\times J_1 \times \dots \times J_N}$ and $\mathcal{Y} \define (y_{i_1, \dots, i_M, j_1, \dots,j_N}) \in \mathbb{C}^{I_1 \times \dots \times I_M\times J_1 \times \dots \times J_N}$, the \emph{Einstein product} of $\mathcal{X}  \star_{N} \mathcal{Y} \in  \mathbb{C}^{I_1 \times \dots \times I_M\times K_1 \times \dots \times K_L}$ is given by 
\begin{eqnarray}\label{eq: Einstein product definition}
(\mathcal{X} \star_{N} \mathcal{Y} )_{i_1, \dots, i_M,k_1,\dots, k_L} \define \sum\limits_{j_1, \dots, j_N} x_{i_1, \dots, i_M, j_1, \dots,j_N}y_{j_1, \dots, j_N, k_1, \dots,k_L}. 
\end{eqnarray}

\begin{definition}\label{def: zero tensor}
A tensor whose entries are all zero is called a \emph{zero tensor}, denoted by $\mathcal{O}$. 
\end{definition}

\begin{definition}\label{def: identity tensor}
An \emph{identity tensor} $\mathcal{I} \in  \mathbb{C}^{I_1 \times \dots \times I_N\times
I_1 \times \dots \times I_N}$ is defined by 
\begin{eqnarray}\label{eq: identity tensor definition}
(\mathcal{I})_{i_1 \times \dots \times i_N\times
j_1 \times \dots \times j_N} \define \prod_{k = 1}^{N} \delta_{i_k, j_k},
\end{eqnarray}
where $\delta_{i_k, j_k} \define 1$ if $i_k  = j_k$; otherwise $\delta_{i_k, j_k} \define 0$.
\end{definition}

In order to define \emph{Hermitian} tensor \cite{ni2019hermitian}, the \emph{conjugate transpose operation} (or \emph{Hermitian adjoint}) of a tensor is specified as follows.  
\begin{definition}\label{def: tensor conjugate transpose}
Given a tensor $\mathcal{X} \define (x_{i_1, \dots, i_M, j_1, \dots,j_N}) \in \mathbb{C}^{I_1 \times \dots \times I_M\times J_1 \times \dots \times J_N}$, its conjugate transpose, denoted by
$\mathcal{X}^{H}$, is defined by
\begin{eqnarray}\label{eq:tensor conjugate transpose definition}
(\mathcal{X}^\mathrm{H})_{ j_1, \dots,j_N,i_1, \dots, i_M}  \define  
 x^*_{i_1, \dots, i_M,j_1, \dots,j_N},
\end{eqnarray}
where the star ``$*$" symbol indicates the complex conjugate of the complex number $x_{i_1, \dots, i_M,j_1, \dots,j_N}$. If a tensor $\mathcal{X}$ satisfies $ \mathcal{X}^\mathrm{H} = \mathcal{X}$, then $\mathcal{X}$ is a \emph{Hermitian tensor}. 
\end{definition}

Following definition is about untiary tensors \cite{ni2019hermitian}.
\begin{definition}\label{def: unitary tensor}
Given a tensor $\mathcal{U} \define (u_{i_1, \dots, i_N, i_1, \dots,i_N}) \in \mathbb{C}^{I_1 \times \dots \times I_N\times I_1 \times \dots \times I_N}$, if
\begin{eqnarray}\label{eq:unitary tensor definition}
\mathcal{U}^\mathrm{H} \star_N \mathcal{U} = \mathcal{U} \star_N \mathcal{U}^\mathrm{H} = \mathcal{I} \in \mathbb{C}^{I_1 \times \dots \times I_N\times I_1 \times \dots \times I_N},
\end{eqnarray}
then $\mathcal{U}$ is a \emph{unitary tensor}. 
\end{definition}

\begin{definition}\label{def: inverse of a tensor}
Given a \emph{square tensor} $\mathcal{X} \define (x_{i_1, \dots, i_N, j_1, \dots,j_N}) \in \mathbb{C}^{I_1 \times \dots \times I_N\times I_1 \times \dots \times I_N}$, if there exists $\mathcal{Y} \in \mathbb{C}^{I_1 \times \dots \times I_N\times I_1 \times \dots \times I_N}$ such that 
\begin{eqnarray}\label{eq:tensor invertible definition}
\mathcal{X} \star_N \mathcal{Y} = \mathcal{Y} \star_N \mathcal{X} = \mathcal{I},
\end{eqnarray}
then $\mathcal{Y}$ is the \emph{inverse} of $\mathcal{X}$. We usually write $\mathcal{Y} \define \mathcal{X}^{-1}$ thereby. 
\end{definition}

We also list other crucial tensor operations here. The \emph{trace} of a square tensor is equivalent to the summation of all diagonal entries such that 
\begin{eqnarray}\label{eq: tensor trace def}
\mathrm{Tr}(\mathcal{X}) \define \sum\limits_{1 \leq i_j \leq I_j,\hspace{0.05cm}j \in [N]} \mathcal{X}_{i_1, \dots, i_N,i_1, \dots, i_N},
\end{eqnarray}
where $[N] \define \{1,2,\cdots,N\}$. The \emph{inner product} of two tensors $\mathcal{X}$, $\mathcal{Y} \in \mathbb{C}^{I_1 \times \dots \times I_N\times J_1 \times \dots \times J_N}$ is given by 
\begin{eqnarray}\label{eq: tensor inner product def}
\langle \mathcal{X}, \mathcal{Y} \rangle \define \mathrm{Tr}\left(\mathcal{X}^H \star_M \mathcal{Y}\right).
\end{eqnarray}

As the matrix eigen-decomposition theorem is crucial in various linear algebra theory and applications, we will have a parallel decomposition theorem for Hermitian tensors. From Theorem 5.2 in~\cite{ni2019hermitian}, every Hermitian tensor $\mathcal{H} \in  \mathbb{C}^{I_1 \times \dots \times I_N \times I_1 \times \dots \times I_N}$ has the following decomposition:
\begin{eqnarray}\label{eq:Hermitian Eigen Decom}
\mathcal{H} &=& \sum\limits_{i=1}^r \lambda_i \mathcal{U}_i  \star_1 \mathcal{U}^{H}_i, \mbox{
~with~~$\langle \mathcal{U}_i, \mathcal{U}_i \rangle =1$ and $\langle \mathcal{U}_i, \mathcal{U}_j \rangle = 0$ for $i \neq j$,}
\end{eqnarray}
where $\lambda_i \in \mathbb{R}$ and $\mathcal{U}_i \in \mathbb{C}^{I_1 \times \dots \times I_N \times 1}$. Here tensors $\mathcal{U}_i$ are orthogonal tensors each other since $\langle \mathcal{U}_i, \mathcal{U}_j \rangle = 0$ for $i \neq j$. The values $\lambda_i$ are named as \emph{Hermitian eigenvalues}, and the minimum integer of $r$ to decompose a Hermitian tensor as in Eq.~\eqref{eq:Hermitian Eigen Decom} is called \emph{Hermitian tensor rank}. In this work, we assume that all Hermitian tensors discussed in this work are full rank, i.e., $r=\prod\limits_{j=1}^N I_j$. A \emph{positive definite} (PD) tensor is a Hermitian tensor with all \emph{Hermitian eigenvalues} are positive. A \emph{semipositive definite} (SPD) tensor is a Hermitian tensor of which all \emph{Hermitian eigenvalues} are nonnegative.  

We use $\preceq$ to repsent  L\"owner ordering between two tensors $\mathcal{A}$ and $\mathcal{B}$ as $\mathcal{A} \preceq \mathcal{B}$, which indicates that $\mathcal{B} - \mathcal{A}$ is a SPD tensor. On the other hand, we use $\succeq$ to repsent  L\"owner ordering between two tensors $\mathcal{A}$ and $\mathcal{B}$ as $\mathcal{A} \succeq \mathcal{B}$, which indicates that $\mathcal{A} - \mathcal{B}$ is a SPD tensor.

\section{Tensor Markov's/Chebyshev’s Inequalities}\label{sec:Tensor Markov Chebyshev’s Inequalities}

Given two Hermitian tensors $\mathcal{A} \in \mathbb{C}^{I_1 \times \cdots \times I_N \times I_1 \times \cdots \times I_N}$ and $\mathcal{B} \in \mathbb{C}^{I_1 \times \cdots \times I_N \times I_1 \times \cdots \times I_N}$, we use $\mathcal{A} \preceq \mathcal{B}$ indicates that the tensor $\mathcal{B} - \mathcal{A}$ is SPD. Such ordering relation among Hermitian tensors can be considered as Loewner ordering. On the other hand, we use  $\mathcal{A} \npreceq \mathcal{B}$ to represent that $\mathcal{B} - \mathcal{A}$ is \emph{not} SPD. This means that $\mathcal{B} - \mathcal{A}$ must contain some negative eigenvalues. 

We will present the following theorem: Markov's inequality for tensors. 
\begin{theorem}[Markov's Inequality for Tensors]\label{thm:Markov's Inequality for Tensors}
Let $\mathcal{A} \in \mathbb{R}^{I_1 \times \cdots \times I_N \times I_1 \times \cdots \times I_N} $ be a PD deterministic tensor, and let $\mathcal{X}  \in \mathbb{R}^{I_1 \times \cdots \times I_N \times I_1 \times \cdots \times I_N}$ be a random Hermitian tensor such that $\mathcal{X} \succeq \mathcal{O}$ almost surely. We have
\begin{eqnarray}\label{eq1:thm:Markov's Inequality for Tensors}
\mathrm{Pr}\left(\mathcal{X} \npreceq \mathcal{A}\right) \leq \mathrm{Tr}\left(\mathbb{E}[\mathcal{X}]\star_N \mathcal{A}^{-1}\right).
\end{eqnarray}
\end{theorem}
\textbf{Proof:}
We have to show the following inequality first. 
\begin{eqnarray}\label{eq2:thm:Markov's Inequality for Tensors}
\mathbbm{1}_{\mathcal{X} \npreceq \mathcal{A}} \leq \mathrm{Tr}\left(\mathcal{A}^{-1/2}\star_N \mathcal{X} \star_N \mathcal{A}^{-1/2}\right),
\end{eqnarray}
where $\mathbbm{1}_{\mathcal{X} \npreceq \mathcal{A}}$ is the indicator function with respect to the condition $\mathcal{X} \npreceq \mathcal{A}$. For the case that $\mathcal{X} \preceq \mathcal{A}$, this inequality provided by Eq.~\eqref{eq2:thm:Markov's Inequality for Tensors} is valid since $\mathbbm{1}_{\mathcal{X} \npreceq \mathcal{A}}$ is equal to zero and the tensor $\mathcal{A}^{-1/2}\star_N \mathcal{X} \star_N \mathcal{A}^{-1/2}$ is a PD tensor.

For the case that $\mathcal{X} \npreceq \mathcal{A}$, we know that at least one tensor $\mathcal{U} \in  \mathbb{C}^{I_1 \times \cdots \times I_N \times 1}$ such that 
\begin{eqnarray}\label{eq3:thm:Markov's Inequality for Tensors}
\mathcal{U}^{\mathrm{T}}\left(\mathcal{A} - \mathcal{X}\right)\mathcal{U} < 0.
\end{eqnarray}
Because the tensor $\mathcal{A}$ is an invertible tensor, we can transform the tensor $\mathcal{U}$ to the tensor $\mathcal{V}$ by setting $\mathcal{V}$$=$$\mathcal{A}^{1/2}\star_N\mathcal{U}$ and obtain the following inequality from Eq.~\eqref{eq3:thm:Markov's Inequality for Tensors}:
\begin{eqnarray}
\mathcal{V}^{\mathrm{T}}\star_N \mathcal{V} - 
\mathcal{V}^{\mathrm{T}}\star_N \mathcal{A}^{-1/2} \star \mathcal{X} \star \mathcal{A}^{1/2} \star_N \mathcal{V} < 0;
\end{eqnarray}
and this is equivalent to 
\begin{eqnarray}\label{eq4:thm:Markov's Inequality for Tensors}
\frac{\mathcal{V}^{\mathrm{T}}\star_N \mathcal{A}^{-1/2} \star \mathcal{X} \star \mathcal{A}^{1/2} \star_N \mathcal{V}}{\mathcal{V}^{\mathrm{T}}\star_N \mathcal{V}}>1.
\end{eqnarray}
By the Rayleigh quotient form of the maximum eigenvalue of a Hermitian tensor, we have 
\begin{eqnarray}
\lambda_{\max}\left(\mathcal{A}^{-1/2} \star \mathcal{X} \star \mathcal{A}^{1/2}\right)&=&\sup\limits_{\mathcal{W} \neq \mathcal{O}}\frac{\mathcal{W}^{\mathrm{T}}\star_N \mathcal{A}^{-1/2} \star \mathcal{X} \star \mathcal{A}^{1/2} \star_N \mathcal{W}}{\mathcal{W}^{\mathrm{T}}\star_N \mathcal{W}} \nonumber\\
&\geq&\frac{\mathcal{V}^{\mathrm{T}}\star_N \mathcal{A}^{-1/2} \star \mathcal{X} \star \mathcal{A}^{1/2} \star_N \mathcal{V}}{\mathcal{V}^{\mathrm{T}}\star_N \mathcal{V}} \nonumber \\
&>&1
\end{eqnarray}
where $\lambda_{\max}\left(\mathcal{A}^{-1/2} \star \mathcal{X} \star \mathcal{A}^{1/2}\right)$ represents the largest eigenvalue for the tensor $\mathcal{A}^{-1/2} \star \mathcal{X} \star \mathcal{A}^{1/2}$. Because $\mathcal{A}^{-1/2} \star \mathcal{X} \star \mathcal{A}^{1/2}$ is a PD tensor, we have
\begin{eqnarray}
\mathrm{Tr}\left(\mathcal{A}^{-1/2} \star \mathcal{X} \star \mathcal{A}^{1/2}\right) \geq \lambda_{\max}\left(\mathcal{A}^{-1/2} \star \mathcal{X} \star \mathcal{A}^{1/2}\right)  > 1.
\end{eqnarray}
This proves Eq.~\eqref{eq2:thm:Markov's Inequality for Tensors}.

By taking the expectation for both sides of Eq.~\eqref{eq2:thm:Markov's Inequality for Tensors}, we have 
\begin{eqnarray}
\mathrm{Pr}\left(\mathcal{X} \npreceq \mathcal{A}\right) &=& \mathbb{E}\mathbbm{1}_{\mathcal{X} \npreceq \mathcal{A}} \nonumber \\
&\leq&  \mathbb{E}\left[\mathrm{Tr}\left(\mathcal{A}^{-1/2}\star_N \mathcal{X} \star_N \mathcal{A}^{-1/2}\right)\right]  \nonumber \\
&=_1& \mathrm{Tr}\left(\mathcal{A}^{-1/2}\star_N \mathbb{E}\left[\mathcal{X}\right] \star_N \mathcal{A}^{-1/2}\right)  \nonumber \\
&=_2& \mathrm{Tr}\left(\mathbb{E}\left[\mathcal{X}\right] \star_N \mathcal{A}^{-1}\right),
\end{eqnarray}
where $=_1$ holds by the linearity of trace and expectation, and  $=_2$ comes from the cyclic multiplication invariance property of trace.
$\hfill \Box$

Before proving tensor Chebyshev's inequality, we need to have the following Lemma.

\begin{lemma}\label{lma:Cheby Lemma}
Let $\mathcal{A}, \mathcal{B} \in \mathbb{C}^{I_1 \times \cdots \times I_N \times I_1 \times \cdots \times I_N}$ be two Hermitian tensors with $\mathcal{A}^2 $$\preceq$$ \mathcal{B}^2$. Then, we have $\left\vert\mathcal{A}\right\vert$$\preceq$$ \left\vert\mathcal{B}\right\vert$. 
\end{lemma}
\textbf{Proof:}
We consider the first case that the tensor $\mathcal{A}$ is an SPD tensor and $\mathcal{B}$ is a PD tensor . From $\mathcal{A}^2 $$\preceq$$ \mathcal{B}^2$, we have
\begin{eqnarray}\label{eq1:lma:Cheby Lemma}
\mathcal{B}^{-1} \star_N \mathcal{A}^{2} \star_N \mathcal{B}^{-1} \preceq \mathcal{I}.
\end{eqnarray}
Then, we have the following inequalities:
\begin{eqnarray}\label{eq2:lma:Cheby Lemma}
1 &\geq_1& \lambda_{\max}\left( \mathcal{B}^{-1} \star_N \mathcal{A}^{2} \star_N \mathcal{B}^{-1}\right) 
= \sigma^2_{\max}\left( \mathcal{B}^{-1} \star_N \mathcal{A}\right) \nonumber \\
&\geq_2& \lambda^2_{\max}\left( \mathcal{B}^{-1} \star_N \mathcal{A}\right) \nonumber \\
&=_3&  \lambda^2_{\max}\left( \mathcal{B}^{-1/2} \star_N \mathcal{A} \star_N \mathcal{B}^{-1/2}\right), 
\end{eqnarray}
where $\geq_1$ comes from Eq.~\eqref{eq1:lma:Cheby Lemma}, $\geq_2$ comes from the Weyl's inequality between eigenvalues and singular values, and $=_3$ comes from the eigenvalue set invariance under the cyclic permutation of tensors. From Eq.~\eqref{eq2:lma:Cheby Lemma}, we have
\begin{eqnarray}\label{eq3:lma:Cheby Lemma}
1 \geq \lambda_{\max}\left( \mathcal{B}^{-1/2} \star_N \mathcal{A}\star_N \mathcal{B}^{-1/2}\right)
\Longleftrightarrow \mathcal{A} \preceq \mathcal{B}.
\end{eqnarray}

Next, we consider the situation that the tensor $\mathcal{A}$ is SPD tensor and $\mathcal{B}$ is a SPD tensor.   From eigen-decomposition of tensor $\mathcal{B}$, we can express $\mathcal{B}$ as $\mathcal{B} $$=$$ \mathcal{U}\star_N \mathcal{D} \star_N \mathcal{U}^{\mathrm{T}}$ , where $\mathcal{U}$ is the unitary tensor and $\mathcal{D}$ is the diagonal tensor. For any $\epsilon >0$, we have $\mathcal{B}+ \epsilon \mathcal{I} $$=$$ \mathcal{U}\star_N (\mathcal{D}+\epsilon \mathcal{I}) \star_N \mathcal{U}^{\mathrm{T}}$, which is PD. From the following relation, 
\begin{eqnarray}
\mathcal{B}^2 &=& \mathcal{U}\star_N \mathcal{D}^2 \star_N \mathcal{U}^{\mathrm{T}} \nonumber \\
&\preceq& \mathcal{U}\star_N \left(\mathcal{D} + \epsilon\mathcal{I}\right)^2 \star_N \mathcal{U}^{\mathrm{T}}  \nonumber \\
&=& \left(\mathcal{B}+ \epsilon \mathcal{I}\right)^2,
\end{eqnarray}
and $\mathcal{A}^2 $$\preceq$$ \mathcal{B}^2$, we have $\mathcal{A}^2 $$\preceq$$ \left(\mathcal{B}+ \epsilon \mathcal{I}\right)^2$. Since $\left(\mathcal{B}+ \epsilon \mathcal{I}\right)$ is a PD tensor, this implies that $\mathcal{A}  $$\preceq$$ \mathcal{B}+ \epsilon \mathcal{I}$. By taking $\epsilon \rightarrow 0$, we have  $\mathcal{A}  $$\preceq$$ \mathcal{B}$.

The remaining case is that both $\mathcal{A}$ and $\mathcal{B}$ are not SPD. Since both $\mathcal{A}$ and $\mathcal{B}$ are Hermitian tensors, we have 
\begin{eqnarray}
\mathcal{A}^2 \preceq \mathcal{B}^2 \Longleftrightarrow 
\left\vert \mathcal{A} \right\vert^2 \preceq \left\vert \mathcal{B} \right\vert^2. 
\end{eqnarray}
Because both $\left\vert \mathcal{A} \right\vert$ and $\left\vert \mathcal{B} \right\vert$ are SPD tensors, we have $\left\vert \mathcal{A} \right\vert $$\preceq$$ \left\vert \mathcal{B} \right\vert$. 
$\hfill \Box$

From Lemma~\ref{lma:Cheby Lemma}, we can have following tensor Chebyshev's inequality.
\begin{theorem}[Chebyshev's Inequality for Tensors]\label{thm:Chebyshev's inequality}
Let $\mathcal{A} \in \mathbb{C}^{I_1 \times \cdots \times I_N \times I_1 \times \cdots \times I_N}$ be a PD deterministic tensor, and let $\mathcal{X} \in \mathbb{C}^{I_1 \times \cdots \times I_N \times I_1 \times \cdots \times I_N}$ be a random Hermitian tensor. Then, we have
\begin{eqnarray}\label{eq1:thm:Chebyshev's inequality}
\mathrm{Pr}\left(\mathcal{X} \npreceq \mathcal{A}\right) \leq \mathrm{Tr}\left(\mathbb{E}[\mathcal{X}^2]\star_N \mathcal{A}^{-2}\right).
\end{eqnarray}
\end{theorem}
\textbf{Proof:}
From Lemma~\ref{lma:Cheby Lemma}, we have $\mathcal{X}^2 $$\preceq$$ \mathcal{A}^2 \Longrightarrow \left\vert \mathcal{X} \right\vert $$\preceq$$ \left\vert \mathcal{A} \right\vert$  This is equivalent that $\left\vert \mathcal{X} \right\vert $$\npreceq$$ \left\vert \mathcal{A} \right\vert\Longrightarrow \mathcal{X}^2 $$\npreceq$$ \mathcal{A}^2$. 

By monotonicity of probability, we have 
\begin{eqnarray}\label{eq2:thm:Chebyshev's inequality}
\mathrm{Pr}\left(\left\vert\mathcal{X}\right\vert \npreceq \mathcal{A}\right)&\leq&
\mathrm{Pr}\left(\mathcal{X}^2 \npreceq \mathcal{A}^2 \right) \nonumber \\
&\leq& \mathrm{Tr}\left(\mathbb{E}\left[\mathcal{X}^2\right]\star_N \mathcal{A}^{-2}\right),
\end{eqnarray}
where the last inequality comes from Theorem~\ref{thm:Markov's Inequality for Tensors}.
$\hfill \Box$

% Lemma Thm 2.3

Actually, we can have more general power, instead of 2, in tensor Chebyshev's inequality. We need the following Lemma~\ref{lma:2.3}.

\begin{lemma}\label{lma:2.3}
Let $\mathcal{A}, \mathcal{B} \in \mathbb{C}^{I_1 \times \cdots \times I_N \times I_1 \times \cdots \times I_N}$ be two SPD tensors with $\mathcal{A} \succeq \mathcal{B}$,  and let $q $$\in$$ [0,1]$. Then, we have
\begin{eqnarray}\label{eq1:lma:2.3}
\mathcal{A}^q \succeq \mathcal{B}^q. 
\end{eqnarray}
\end{lemma}
\textbf{Proof:}
Let us define the following set, $\mathrm{S}$, for the real number $q$:
\begin{eqnarray}\label{eq2:lma:2.3}
\mathrm{S} &=& \{q \in \mathbb{R}: \mathcal{A}^q \succeq\mathcal{B}^q\}.
\end{eqnarray}
It is clear that both $0$ and $1$ belong to the set $\mathrm{S}$. We wish to show that $[0,1] \in \mathrm{S}$. Since dyadic numbers are dense in $[0,1]$, this is equivalent to show the following:
\begin{eqnarray}\label{eq3:lma:2.3}
q,r \in \mathrm{S} \Longrightarrow \frac{q+r}{2} \in \mathrm{S}.
\end{eqnarray}

If $q \in \mathrm{S}$, we have 
\begin{eqnarray}
\mathcal{A}^q \succeq \mathcal{B}^q \Longrightarrow \mathcal{A}^{-q/2} \star_N \mathcal{B}^q 
\star_N \mathcal{A}^{-q/2} \preceq \mathcal{I}.
\end{eqnarray}
Then, we can have the following 
\begin{eqnarray}\label{eq4:lma:2.3}
\lambda^2_{\max}\left(\mathcal{B}^{q/2}\star_N \mathcal{A}^{-q/2}\right)&=&
\lambda_{\max}\left(\left(\mathcal{B}^{q/2}\star_N \mathcal{A}^{-q/2}\right)^{\mathrm{H}}\star_N \left(\mathcal{B}^{q/2}\star_N \mathcal{A}^{-q/2}\right)\right) \nonumber \\
&\leq& \lambda_{\max}\left(\mathcal{A}^{-q/2}\star_N\mathcal{B}^{q}\star_N \mathcal{A}^{-q/2}\right) \leq 1
\end{eqnarray}
This implies that $\lambda_{\max}\left(\mathcal{B}^{q/2}\star_N \mathcal{A}^{-q/2}\right) \leq 1$. Because we also have $r \in \mathrm{S}$, similarly, we have $\lambda_{\max}\left(\mathcal{B}^{r/2}\star_N \mathcal{A}^{-r/2}\right) \leq 1$.

Therefore, we have 
\begin{eqnarray}\label{eq5:lma:2.3}
1 &\geq& \lambda_{\max}\left(\left(\mathcal{B}^{q/2}\star_N \mathcal{A}^{-q/2}\right)^{\mathrm{H}}\star_N \left(\mathcal{B}^{r/2}\star_N \mathcal{A}^{-r/2}\right)\right)\nonumber \\
&=& \lambda_{\max}\left(\mathcal{A}^{-q/2}\star_N \mathcal{B}^{(q+r)/2}\star_N \mathcal{A}^{-r/2}\right) \nonumber \\
&\geq& \lambda_{\max}\left(\mathcal{A}^{-(q+r)/2}\star_N \mathcal{B}^{(q+r)/2}\star_N \mathcal{A}^{-(q+r)/2}\right),
\end{eqnarray}
which implies that $\mathcal{A}^{(q+r)/2} $$\succeq$$ \mathcal{B}^{(q+r)/2}$. This Lemma is proved since $(q+r)/2 \in \mathrm{S}$.
$\hfill \Box$

With this Lemma~\ref{lma:2.3}, we can have the following generalized Chebyshev's Inequality for Tensors with general exponent. 
\begin{theorem}[Generalized Chebyshev's Inequality for Tensors]\label{thm:Chebyshev's inequality q}
Let $\mathcal{A} $$\in$$ \mathbb{C}^{I_1 \times \cdots \times I_N \times I_1 \times \cdots \times I_N}$ be a PD deterministic tensor, and let $\mathcal{X} \in \mathbb{C}^{I_1 \times \cdots \times I_N \times I_1 \times \cdots \times I_N}$ be a random Hermitian tensor. Given $p \geq 1$, then we have
\begin{eqnarray}\label{eq1:thm:Chebyshev's inequality q}
\mathrm{Pr}\left(\mathcal{X} \npreceq \mathcal{A}\right) \leq \mathrm{Tr}\left(\mathbb{E}[\left\vert\mathcal{X}\right\vert^p]\star_N \mathcal{A}^{-p}\right).
\end{eqnarray}
\end{theorem}
\textbf{Proof:}
We have $\left\vert \mathcal{X} \right\vert^p \preceq \mathcal{A}^p \Longrightarrow \left\vert \mathcal{X} \right\vert \preceq \mathcal{A}$ by setting $q=1/p$ in  Lemma~\ref{lma:2.3}. Again, by monotonicity of probability and Theorem~\ref{thm:Markov's Inequality for Tensors}, we have
\begin{eqnarray}
\mathrm{Pr}\left(\left\vert \mathcal{X}\right\vert \npreceq \mathcal{A} \right) \leq 
\mathrm{Pr}\left(\left\vert \mathcal{X}\right\vert^p \npreceq \mathcal{A}^p \right) \leq 
\mathrm{Tr}\left(\mathbb{E}\left[\left\vert \mathcal{X}\right\vert^p \right]\star_N \mathcal{A}^p \right).
\end{eqnarray}
$\hfill \Box$

Following Lemma is used to associate Loewner ordering with Theorem~\ref{thm:Chebyshev's inequality q}.
\begin{lemma}\label{lma:Loewner ordering with Markov Cheb inequalities}
Given the following random PD tensors $\mathcal{X}, \mathcal{Y}, \mathcal{Z} $$\in$$ \mathbb{C}^{I_1 \times \cdots \times I_N \times I_1 \times \cdots \times I_N}$ with the relation $\mathcal{X} $$\preceq$$ \mathcal{Y} $$\preceq$$ \mathcal{Z}$, and a deterministic PD tensor $\mathcal{C}$, we have 
\begin{eqnarray}\label{eq1:lma:Loewner ordering with Markov Cheb inequalities}
\mathrm{Pr}\left(\mathcal{Y} \npreceq \mathcal{C}\right) &\leq& \mathrm{Tr}\left(\mathbb{E}\left[\mathcal{Z}^q\right] \star_N \mathcal{C}^{-1}\right),
\end{eqnarray}
where $q \geq 1$. We also have 
\begin{eqnarray}\label{eq2:lma:Loewner ordering with Markov Cheb inequalities}
\mathrm{Pr}\left(\mathcal{X} \npreceq \mathcal{C}\right) &\leq& \mathrm{Tr}\left(\mathbb{E}\left[\mathcal{Y}^q\right] \star_N \mathcal{C}^{-1}\right).
\end{eqnarray}
\end{lemma}
\textbf{Proof:}
From Theorem~\ref{thm:Chebyshev's inequality q}, we have
\begin{eqnarray}\label{eq3:lma:Loewner ordering with Markov Cheb inequalities}
\mathrm{Pr}\left(\mathcal{Y} \npreceq \mathcal{C}\right) &\leq& \mathrm{Tr}\left(\mathbb{E}\left[\mathcal{Y}^q\right] \star_N \mathcal{C}^{-q}\right) \nonumber \\
&\leq&  \mathrm{Tr}\left(\mathbb{E}\left[\mathcal{Z}^q\right] \star_N \mathcal{C}^{-q}\right)
\end{eqnarray}
where the last inequality comes from the tensor monotone property of the function $x^{q}$ for $q \geq 1$. Then, we have Eq.~\eqref{eq1:lma:Loewner ordering with Markov Cheb inequalities}.

Since $\mathcal{X} $$\preceq$$ \mathcal{Y}$, by probability monotone property, we have 
\begin{eqnarray}\label{eq4:lma:Loewner ordering with Markov Cheb inequalities}
\mathrm{Pr}\left(\mathcal{X} \npreceq \mathcal{C}\right) &\leq& \mathrm{Pr}\left(\mathcal{Y} \npreceq \mathcal{C}\right) \nonumber \\
&\leq&  \mathrm{Tr}\left(\mathbb{E}\left[\mathcal{Y}^q\right] \star_N \mathcal{C}^{-q}\right)
\end{eqnarray}
where the last inequality comes from Theorem~\ref{thm:Chebyshev's inequality q}.
$\hfill \Box$

\section{Tail Bounds for Bivariate Random Tensor Means Based on L\"oewner Ordering}\label{sec:Tail Bounds for Bivariate Random Tensors Mean Based on Loewner Ordering} 

In this section, we will establish several tail bounds for bivariate random tensor means based on Loewner ordering. We will introduce the basic notion about \emph{tensor monotone} and its properties. 

A real continuous function $f$ defined on $(0, \infty)$ will be named as a \emph{tensor monotone increasing} function if we have
\begin{eqnarray}\label{eq:tensor monotone def}
\mathcal{A} \succeq \mathcal{B} \succ \mathcal{O} \Longrightarrow 
f(\mathcal{A}) \succeq f(\mathcal{B}),
\end{eqnarray}
where $\mathcal{A}, \mathcal{B}$ are Hermitian tensors. Similarly, a function $f$ is named as a \emph{tensor monotone decreasing} function if $-f$ is a \emph{tensor monotone increasing} function. Besides, a real continuous function $f$ defined on $(0, \infty)$ will be named as a \emph{tensor convex} function if we have
\begin{eqnarray}\label{eq:tensor convex def}
\lambda f(\mathcal{A}) + (1 - \lambda)f(\mathcal{B}) \succeq f(\lambda \mathcal{A} + (1 - \lambda) \mathcal{B}),
\end{eqnarray}
where $\mathcal{A}, \mathcal{B}$ are PD tensors and $\lambda \in [0,1]$. Then, we will define three sets of positive functions as:
\begin{eqnarray}\label{eq:three sets of functions}
\mbox{TMI}&=&\{f: \mbox{tensor monotone increasing on $(0, \infty)$, $f  $$>$$ 0$}\};\nonumber \\
\mbox{TMD}&=&\{g: \mbox{tensor monotone decreasing on $(0, \infty)$, $g $$>$$ 0$}\};\nonumber \\
\mbox{TC}&=&\{h: \mbox{tensor convex on $(0, \infty)$, $h $$>$$ 0$}\}.
\end{eqnarray}
Moreover, we define the following three sets of positive functions as:
\begin{eqnarray}\label{eq:three sets of functions with one}
\mbox{TMI}^{1}&=&\{f: \mbox{tensor monotone increasing on $(0, \infty)$, $f $$>$$ 0$ and $f(1) $$=$$ 1$}\};\nonumber \\
\mbox{TMD}^{1}&=&\{g: \mbox{tensor monotone decreasing on $(0, \infty)$, $g $$>$$ 0$ and $g(1) $$=$$ 1$}\};\nonumber \\
\mbox{TC}^{1}&=&\{h: \mbox{tensor convex on $(0, \infty)$, $h $$>$$ 0$ and $h(1) $$=$$ 1$}\}.
\end{eqnarray}

\subsection{Connection Functions Come From $\mbox{TMI}^{1}$ or $\mbox{TMD}^{1}$}\label{sec:Connection Functions are TCI or TMD}

The bivariate PD tensor means for tensors $\mathcal{X} $$\in$$ \mathbb{C}^{I_1 \times \cdots \times I_N \times I_1 \times \cdots \times I_N}$ and $\mathcal{Y} $$\in$$ \mathbb{C}^{I_1 \times \cdots \times I_N \times I_1 \times \cdots \times I_N}$ with respect to the function $f \in \mbox{TMI}^{1}$, denoted by $\mathcal{X} \#_f \mathcal{Y}$, is defined as~\footnote{Here, we adopt Kubo-Ando's sense operator mean for tensors mean~\cite{kubo1980means}.} 
\begin{eqnarray}
\mathcal{X} \#_f \mathcal{Y} &\define& \mathcal{X}^{1/2} \star_N f\left(\mathcal{X}^{-1/2} \star_N \mathcal{Y} \star_N \mathcal{X}^{-1/2}\right)\star_N \mathcal{X}^{1/2}.
\end{eqnarray}

Following theorem is about tail bounds for random tensors $\mathcal{X} $$\#_f$$ \mathcal{Y}$ and $\mathcal{X}^q $$\#_f$$ \mathcal{Y}^q$ with respect to exponent $q$. Let us define the generalized product operation, denoted by $\acute{\prod}_{k=1}^n$, when the index upper bound is less than the index lower bound:
\begin{eqnarray}\label{eq4-1-2:thm:3.2}
\acute{\prod}_{k=1}^n a_i \define \begin{cases}
\prod_{k=1}^n a_i, ~\mbox{if $n \geq 1$};        \\
        1,~\mbox{if $n = 0$}.
        \end{cases}
\end{eqnarray}
where $a_i$ is the $i$-th real number.

The function $f$ has \emph{power monotone increasing} (pmi for abbreviation) property if it satisfies the following:
\begin{eqnarray}\label{eq:pmi def}
f^q(x) \leq f(x^q),
\end{eqnarray}
where $x  $$>$$ 0$ and $q $$\geq$$ 1$. On the other hand, the function $f$ has \emph{power monotone decreasing} (pmd for abbreviation) property if it satisfies the following:
\begin{eqnarray}\label{eq:pmd def}
f^q(x) \geq f(x^q). 
\end{eqnarray}

\begin{theorem}\label{thm:3.2}
Given two random PD tensors $\mathcal{X} $$\in$$ \mathbb{C}^{I_1 \times \cdots \times I_N \times I_1 \times \cdots \times I_N}$, $\mathcal{Y} $$\in$$ \mathbb{C}^{I_1 \times \cdots \times I_N \times I_1 \times \cdots \times I_N}$ and a PD determinstic tensor $\mathcal{C}$, if $q=2^n q_0 \geq 1$ with $1 \leq q_0 \leq 2$ and $n \in \mathbb{N}$, we set $\mathcal{Z}_{k-1} \define \mathcal{X}^{-2^{k-2}} \mathcal{Y}^{2^{k-1}}\mathcal{X}^{-2^{k-2}}$ for $k=1,2,\cdots,n$. We assume that $\mathcal{X} $$\#_f$$ \mathcal{Y} $$\succeq$$ \mathcal{I}$ almost surely with $f $$\in$$\mbox{TMI}^{1}$, we have 
\begin{eqnarray}\label{eq1-1:thm:3.2}
\mathrm{Pr}\left(\mathcal{X}^q \#_f \mathcal{Y}^q \npreceq \mathcal{C} \right)
&\leq& \mathrm{Tr}\left(\mathbb{E}\left[\left(\Psi_{upper}\left(q,f,\mathcal{X},\mathcal{Y}\right)\lambda^{q-1}_{\min}\left(\mathcal{X}\#_f\mathcal{Y}\right)\mathcal{X}\#_f\mathcal{Y}\right)^p\right]\star_N \mathcal{C}^{-1} \right),
\end{eqnarray}
and
\begin{eqnarray}\label{eq1-2:thm:3.2}
\mathrm{Pr}\left(\Psi_{lower}\left(q,f,\mathcal{X},\mathcal{Y}\right)\lambda^{q-1}_{\max}\left(\mathcal{X}\#_f\mathcal{Y}\right)\mathcal{X}\#_f\mathcal{Y} \npreceq \mathcal{C} \right)
&\leq& \mathrm{Tr}\left(\mathbb{E}\left[\left(\mathcal{X}^q \#_f \mathcal{Y}^q\right)^p\right]\star_N \mathcal{C}^{-1} \right),
\end{eqnarray}
where $\Psi_{lower}\left(q,f,\mathcal{X},\mathcal{Y}\right)$ and $\Psi_{upper}\left(q,f,\mathcal{X},\mathcal{Y}\right)$ are two positive numbers defined by
\begin{eqnarray}
\Psi_{lower}\left(q,f,\mathcal{X},\mathcal{Y}\right)&\define&\lambda_{\min}\left(f^{-q_0}\left(\mathcal{Z}_n\right)f\left(\mathcal{Z}^{q_0}_n\right)\right) \acute{\prod}_{k=1}^n \lambda_{\min}\left(f^{-2}\left(\mathcal{Z}_{k-1}\right)f\left(\mathcal{Z}_{k-1}^{2}\right)\right) \nonumber \\
\Psi_{upper}\left(q,f,\mathcal{X},\mathcal{Y}\right)&\define&\lambda_{\max}\left(f^{-q_0}\left(\mathcal{Z}_n\right)(f\left(\mathcal{Z}_n^{q_0}\right)\right)\acute{\prod}_{k=1}^n \lambda_{\max}\left(f^{-2}\left(\mathcal{Z}_{k-1}\right)f\left(\mathcal{Z}_{k-1}^{2}\right)\right).
\end{eqnarray}
Note that the definition of $\acute{\prod}$ is provided by Eq.~\eqref{eq4-1-2:thm:3.2}.

For $0 < q  \leq1$, we have 
\begin{eqnarray}\label{eq2-1:thm:3.2}
\mathrm{Pr}\left(\mathcal{X}^q \#_f \mathcal{Y}^q \npreceq \mathcal{C} \right)
&\leq& \mathrm{Tr}\left(\mathbb{E}\left[\left(\lambda_{\min}\left(f^{-q}\left(\mathcal{Z}_0\right)(f\left(\mathcal{Z}_0^q\right)\right)\lambda^{q-1}_{\min}\left(\mathcal{X}\#_f\mathcal{Y}\right)\mathcal{X}\#_f\mathcal{Y}\right)^p\right]\star_N \mathcal{C}^{-1} \right),
\end{eqnarray}
and
\begin{eqnarray}\label{eq2-2:thm:3.2}
\mathrm{Pr}\left( \lambda_{\max}\left(f^{-q}\left(\mathcal{Z}_0\right)f\left(\mathcal{Z}_0^q\right)\right)\lambda^{q-1}_{\max}\left(\mathcal{X}\#_f\mathcal{Y}\right)\mathcal{X}\#_f\mathcal{Y} \npreceq \mathcal{C} \right)
&\leq& \mathrm{Tr}\left(\mathbb{E}\left[\left(\mathcal{X}^q \#_f \mathcal{Y}^q\right)^p\right]\star_N \mathcal{C}^{-1} \right).
\end{eqnarray}
where $p \geq 1$. 
\end{theorem}
\textbf{Proof:}
Since all multiplications between tensors are $\star_N$, we will remove these notations in this proof for simplification. We begin with the case for $q \geq 1$. We will separte the region of $q \geq 1$ into $1 \leq q \leq 2$ and $q \geq 2$. 

For the subregion $1 \leq q \leq 2$, $r \define 2-q$, and $\mathcal{X} $$\#_f$$ \mathcal{Y} $$\succeq$$ \mathcal{I}$, we have
\begin{eqnarray}\label{eq4:thm:3.2}
\mathcal{X}^q \#_f \mathcal{Y}^q &=& \mathcal{X}^{\frac{q}{2}}f\left(\mathcal{X}^{\frac{1-q}{2}}\mathcal{Z}_0\mathcal{X}^{\frac{1}{2}}\mathcal{Y}^{-r}\mathcal{X}^{\frac{1}{2}}\mathcal{Z}_0\mathcal{X}^{\frac{1-q}{2}}\right)\mathcal{X}^{\frac{q}{2}} \nonumber \\
&=& \mathcal{X}^{\frac{q}{2}}f\left(\mathcal{X}^{\frac{1-q}{2}}\mathcal{Z}_0\mathcal{X}^{\frac{1}{2}}\left(\mathcal{X}^{\frac{-1}{2}} \mathcal{Z}_0^{-1} \mathcal{X}^{\frac{-1}{2}}\right)^{r}\mathcal{X}^{\frac{1}{2}}\mathcal{Z}_0\mathcal{X}^{\frac{1-q}{2}}\right)\mathcal{X}^{\frac{q}{2}} \nonumber \\
&=& \mathcal{X}^{\frac{1}{2}}\mathcal{X}^{\frac{1-r}{2}}f\left(\mathcal{X}^{\frac{r-1}{2}}\mathcal{Z}_0
\left(\mathcal{X} \#_{x^r} \mathcal{Z}_0^{-1}\right)
\mathcal{Z}_0\mathcal{X}^{\frac{r-1}{2}}\right)\mathcal{X}^{\frac{1-r}{2}}\mathcal{X}^{\frac{1}{2}}\nonumber \\
&=& \mathcal{X}^{\frac{1}{2}}\left[\mathcal{X}^{1-r}\#_f \left(\mathcal{Z}_0\left(\mathcal{X} \#_{x^r} \mathcal{Z}_0^{-1}\right)\mathcal{Z}_0\right)\right]\mathcal{X}^{\frac{1}{2}}\nonumber \\
&\succeq_1&  \mathcal{X}^{\frac{1}{2}}\left[f^{1-q}(\mathcal{Z})\#_f \left(\mathcal{Z}\left(f^{-1}(\mathcal{Z})\#_{x^r} \mathcal{Z}_0^{-1}\right)\mathcal{Z}_0\right)\right]\mathcal{X}^{\frac{1}{2}}\nonumber \\
&=& \mathcal{X}^{\frac{1}{2}}\left[f^{1-q}\left(\mathcal{Z}_0\right)(f\left(\mathcal{Z}_0^q \right)\right]\mathcal{X}^{\frac{1}{2}}\nonumber \\
&\succeq& \lambda_{\min}\left(f^{-q}\left(\mathcal{Z}_0\right)(f\left(\mathcal{Z}_0^q\right)\right)\mathcal{X}\#_f\mathcal{Y},
\end{eqnarray}
where $\succeq_1$ we applies Lemma~\ref{lma:2.3} based on $\mathcal{X} $$\#_f$$ \mathcal{Y} $$\succeq$$ \mathcal{I} \Longleftrightarrow \mathcal{X} \succeq f^{-1}(\mathcal{Z}_0)$. If $q \geq 2$, by finding some natural number $n$ such that $q=2^n q_0$ with $1 \leq q_0 \leq 2$ and iterating the relation provided by Eq.~\eqref{eq4:thm:3.2} dyadically with respect to $q$, we have the relation:
\begin{eqnarray}\label{eq4-1:thm:3.2}
\mathcal{X}^q \#_f \mathcal{Y}^q &\succeq& \lambda_{\min}\left(f^{-q_0}\left(\mathcal{Z}_n\right)f\left(\mathcal{Z}^{q_0}_n\right)\right) \prod_{k=1}^n \lambda_{\min}\left(f^{-2}\left(\mathcal{Z}_{k-1}\right)f\left(\mathcal{Z}_{k-1}^{2}\right)\right) \mathcal{X} \#_f \mathcal{Y}.
\end{eqnarray}
By combining Eq.~\eqref{eq4:thm:3.2} and Eq.~\eqref{eq4-1:thm:3.2}, for $q \geq 1$, we have 
\begin{eqnarray}\label{eq4-1-1:thm:3.2}
\mathcal{X}^q \#_f \mathcal{Y}^q &\succeq& \lambda_{\min}\left(f^{-q_0}\left(\mathcal{Z}_n\right)f\left(\mathcal{Z}^{q_0}_n\right)\right) \acute{\prod}_{k=1}^n \lambda_{\min}\left(f^{-2}\left(\mathcal{Z}_{k-1}\right)f\left(\mathcal{Z}_{k-1}^{2}\right)\right) \mathcal{X} \#_f \mathcal{Y}.
\end{eqnarray}

Given any positive real number, say $\beta$, we can replace $\mathcal{X}$ and $\mathcal{Y}$ in Eq.~\eqref{eq4-1:thm:3.2} with 
$\beta^{-1}\mathcal{X}$ and $\beta^{-1}\mathcal{Y}$ to get
\begin{eqnarray}\label{eq4-2:thm:3.2}
\mathcal{X}^q \#_f \mathcal{Y}^q &\succeq& \lambda_{\min}\left(f^{-q_0}\left(\mathcal{Z}_n\right)f\left(\mathcal{Z}^{q_0}_n\right)\right) \acute{\prod}_{k=1}^n \lambda_{\min}\left(f^{-2}\left(\mathcal{Z}_{k-1}\right)f\left(\mathcal{Z}_{k-1}^{2}\right)\right) \mathcal{X} \#_f \mathcal{Y}.
\end{eqnarray}
We can select $\beta $$=$$ \lambda_{\max}\left(\mathcal{X}\#_f\mathcal{Y}\right)$ to associate $\beta$ with 
$\mathcal{X}\#_f\mathcal{Y}$ and maximize R.H.S. of Eq.~\eqref{eq4-2:thm:3.2} in the sense of L\"oewner ordering. Then, we have 
\begin{eqnarray}\label{eq4-3:thm:3.2}
\mathcal{X}^q \#_f \mathcal{Y}^q\succeq\lambda_{\min}\left(f^{-q_0}\left(\mathcal{Z}_n\right)f\left(\mathcal{Z}^{q_0}_n\right)\right) \acute{\prod}_{k=1}^n \lambda_{\min}\left(f^{-2}\left(\mathcal{Z}_{k-1}\right)f\left(\mathcal{Z}_{k-1}^{2}\right)\right)\lambda^{q-1}_{\max}\left(\mathcal{X}\#_f\mathcal{Y}\right)\mathcal{X}\#_f\mathcal{Y}.
\end{eqnarray}

By replacing $f $$\in$$ \mbox{TMI}^{1}$, $\mathcal{X}$ and $\mathcal{Y}$ in Eq.~\eqref{eq4-3:thm:3.2} with  $f^{\ast} \define f^{-1}(x^{-1})$ for $x $$\in$$ (0,\infty)$, $\mathcal{X}^{-1}$ and $\mathcal{Y}^{-1}$, we have 
\begin{eqnarray}\label{eq5:thm:3.2}
\mathcal{X}^{-q} \#_{f^\ast} \mathcal{Y}^{-q} &\succeq& \lambda_{\min}\left((f^{\ast})^{-q_0}\left(\mathcal{Z}_n'\right)f^{\ast}\left(\mathcal{Z}_n'^{q_0}\right)\right)\acute{\prod}_{k=1}^n \lambda_{\min}\left((f^\ast)^{-2}\left(\mathcal{Z}'_{k-1}\right)(f^{\ast})\left(\mathcal{Z}_{k-1}^{'2}\right)\right)\nonumber \\
& & \cdot \lambda^{q-1}_{\max}\left(\mathcal{X}^{-1}\#_{f^{\ast}}\mathcal{Y}^{-1}\right)\mathcal{X}^{-1}\#_{f^{\ast}}\mathcal{Y}^{-1}.
\end{eqnarray}
where $\mathcal{Z}'_{k-1} \define \mathcal{X}^{2^{k-2}} \mathcal{Y}^{-2^{k-1}}\mathcal{X}^{2^{k-2}}$ for $k=1,2,\cdots,n$. From the definition of $f^{\ast}$, we have 
\begin{eqnarray}\label{eq5-1:thm:3.2}
\mathcal{X}^{-1}\#_{f^{\ast}}\mathcal{Y}^{-1}&=&\left(\mathcal{X}\#_{f}\mathcal{Y}\right)^{-1}, \nonumber \\
\mathcal{X}^{-q}\#_{f^{\ast}}\mathcal{Y}^{-q}&=&\left(\mathcal{X}^q \#_{f} \mathcal{Y}^q\right)^{-1},\nonumber \\
(f^{\ast})^{-q_0}\left(\mathcal{Z}_n'\right)f^{\ast}\left(\mathcal{Z}_n'^{q_0}\right)&=&f^{-1}\left(\mathcal{Z}_n^{q_0}\right)f^{q_0}\left(\mathcal{Z}_n\right),\nonumber \\
(f^{\ast})^{-2}\left(\mathcal{Z}_{k-1}'\right)f^{\ast}\left(\mathcal{Z}_{k-1}'^{2}\right)&=&f^{-1}\left(\mathcal{Z}_{k-1}^{2}\right)f^{2}\left(\mathcal{Z}_{k-1}\right),~\mbox{for $k=1,2,\cdots,n$;}
\end{eqnarray} 
then, by applying Eq.~\eqref{eq5-1:thm:3.2} to Eq.~\eqref{eq5:thm:3.2}, we obtain the following:
\begin{eqnarray}\label{eq5-2:thm:3.2}
\mathcal{X}^q \#_f \mathcal{Y}^q &\preceq& \lambda_{\max}\left(f^{-q_0}\left(\mathcal{Z}_n\right)(f\left(\mathcal{Z}_n^{q_0}\right)\right)\acute{\prod}_{k=1}^n \lambda_{\max}\left(f^{-2}\left(\mathcal{Z}_{k-1}\right)f\left(\mathcal{Z}_{k-1}^{2}\right)\right)\nonumber \\
&& \cdot \lambda^{q-1}_{\min}\left(\mathcal{X}\#_f\mathcal{Y}\right)\mathcal{X}\#_f\mathcal{Y}.
\end{eqnarray} 
By combining Eq.~\eqref{eq4-3:thm:3.2} and Eq.~\eqref{eq5-2:thm:3.2}, for $q \geq 1$, we have 
\begin{eqnarray}\label{eq5-3:thm:3.2}
\Psi_{lower}\left(q,f,\mathcal{X},\mathcal{Y}\right)\lambda^{q-1}_{\max}\left(\mathcal{X}\#_f\mathcal{Y}\right)\mathcal{X}\#_f\mathcal{Y} &\preceq& \mathcal{X}^q \#_f \mathcal{Y}^q  \nonumber \\
&\preceq &\Psi_{upper}\left(q,f,\mathcal{X},\mathcal{Y}\right)\lambda^{q-1}_{\min}\left(\mathcal{X}\#_f\mathcal{Y}\right)\mathcal{X}\#_f\mathcal{Y}.
\end{eqnarray}

By applying Lemma~\ref{lma:Loewner ordering with Markov Cheb inequalities}, we have the desired bounds provided by Eq.~\eqref{eq1-1:thm:3.2} and Eq.~\eqref{eq1-2:thm:3.2}.

Now, we will consider the case for $0 < q \leq 1$. 
\begin{eqnarray}\label{eq6:thm:3.2}
\mathcal{X}^q \#_f \mathcal{Y}^q
&=& \mathcal{X}^{\frac{q}{2}}f\left(\mathcal{X}^{\frac{-q}{2}}\left(\mathcal{X}^{\frac{1}{2}}\mathcal{Z}_0\mathcal{X}^{\frac{1}{2}}\right)^{q}\mathcal{X}^{\frac{-q}{2}}\right)\mathcal{X}^{\frac{q}{2}} \nonumber \\
&=& \mathcal{X}^{\frac{q}{2}}f\left(\mathcal{X}^{\frac{1-q}{2}}\left(\mathcal{X}^{-1}\#_{x^q}\mathcal{Z}_0\right)\mathcal{X}^{\frac{1-q}{2}}\right)\mathcal{X}^{\frac{q}{2}}\nonumber \\
&=& \mathcal{X}^{\frac{1}{2}}\left[\mathcal{X}^{-(1-q)}\#_f \left(\mathcal{X}^{-1}\#_{x^q}\mathcal{Z}_0\right)\right]\mathcal{X}^{\frac{1}{2}}\nonumber \\
&\preceq_1& \mathcal{X}^{\frac{1}{2}}\left[f^{(1-q)}\left(\mathcal{Z}_0\right)\#_f \left(f\left(\mathcal{Z}_0\right)\#_{x^q}\mathcal{Z}_0\right)\right]\mathcal{X}^{\frac{1}{2}}\nonumber \\
&=& \mathcal{X}^{\frac{1}{2}}\left[f^{1-q}\left(\mathcal{Z}_0\right)(f\left(\mathcal{Z}_0^q\right)\right]\mathcal{X}^{\frac{1}{2}}\nonumber \\
&\preceq& \lambda_{\max}\left(f^{-q}\left(\mathcal{Z}_0\right)(f\left(\mathcal{Z}_0^q\right)\right)\mathcal{X}\#_f\mathcal{Y},
\end{eqnarray}
where $\succeq_1$ we applies  Lemma~\ref{lma:2.3} based on $\mathcal{X} $$\#_f$$ \mathcal{Y} $$\succeq$$ \mathcal{I} \Longleftrightarrow \mathcal{X} \succeq f^{-1}(\mathcal{Z})$ with $0 \leq 1-q \leq 1$. Given any positive real number, say $\beta$, we also can replace $\mathcal{X}$ and $\mathcal{Y}$ in Eq.~\eqref{eq6:thm:3.2} with 
$\beta^{-1}\mathcal{X}$ and $\beta^{-1}\mathcal{Y}$ to get
\begin{eqnarray}\label{eq6-1:thm:3.2}
\mathcal{X}^q \#_f \mathcal{Y}^q &\preceq& \lambda_{\min}\left(f^{-q}\left(\mathcal{Z}_0\right)(f\left(\mathcal{Z}_0^q\right)\right)\beta^{q-1}\mathcal{X}\#_f\mathcal{Y}.
\end{eqnarray}
We can select $\beta $$=$$ \lambda_{\min}\left(\mathcal{X}\#_f\mathcal{Y}\right)$ to associate $\beta$ with 
$\mathcal{X}\#_f\mathcal{Y}$ and minimize R.H.S. of Eq.~\eqref{eq6-1:thm:3.2} in the sense of Loewner ordering. Then, we have 
\begin{eqnarray}\label{eq6-2:thm:3.2}
\mathcal{X}^q \#_f \mathcal{Y}^q &\preceq& \lambda_{\min}\left(f^{-q}\left(\mathcal{Z}_0\right)f\left(\mathcal{Z}_0^q\right)\right)\lambda^{q-1}_{\min}\left(\mathcal{X}\#_f\mathcal{Y}\right)\mathcal{X}\#_f\mathcal{Y}.
\end{eqnarray}

By replacing $f $$\in$$ \mbox{TMI}^{1}$, $\mathcal{X}$ and $\mathcal{Y}$ in Eq.~\eqref{eq6-2:thm:3.2} with  $f^{\ast} \define f^{-1}(x^{-1})$ for $x $$\in$$ (0,\infty)$, $\mathcal{X}^{-1}$ and $\mathcal{Y}^{-1}$, we have 
\begin{eqnarray}\label{eq6-3:thm:3.2}
\mathcal{X}^q \#_f \mathcal{Y}^q &\succeq& \lambda_{\max}\left(f^{-q}\left(\mathcal{Z}_0\right)f\left(\mathcal{Z}_0^q\right)\right)\lambda^{q-1}_{\max}\left(\mathcal{X}\#_f\mathcal{Y}\right)\mathcal{X}\#_f\mathcal{Y}.
\end{eqnarray}
Therefore, combining Eq.~\eqref{eq6-2:thm:3.2} and Eq.~\eqref{eq6-3:thm:3.2}, we have 
\begin{eqnarray}\label{eq6-4:thm:3.2}
\lambda_{\max}\left(f^{-q}\left(\mathcal{Z}_0\right)f\left(\mathcal{Z}_0^q\right)\right)\lambda^{q-1}_{\max}\left(\mathcal{X}\#_f\mathcal{Y}\right)\mathcal{X}\#_f\mathcal{Y} &\preceq& \mathcal{X}^q \#_f \mathcal{Y}^q  \nonumber \\
&\preceq &  \lambda_{\min}\left(f^{-q}\left(\mathcal{Z}_0\right)(f\left(\mathcal{Z}_0^q\right)\right)\lambda^{q-1}_{\min}\left(\mathcal{X}\#_f\mathcal{Y}\right)\mathcal{X}\#_f\mathcal{Y},\nonumber \\
\end{eqnarray}
where $0 \leq q \leq 1$. By applying Lemma~\ref{lma:Loewner ordering with Markov Cheb inequalities}, we have the desired bounds provided by Eq.~\eqref{eq2-1:thm:3.2} and Eq.~\eqref{eq2-2:thm:3.2}.
$\hfill \Box$

We can have the following Corollary~\ref{cor:3.3} derived from Theorem~\ref{thm:3.2} with simpler formats by assuming pmi or pmd for the connection function $f$.
\begin{corollary}\label{cor:3.3}
Given same conditions provided by Theorem~\ref{thm:3.2} with the function $f$ satisfying pmi propertity, we have
\begin{eqnarray}\label{eq2-3:cor:3.3}
\mathrm{Pr}\left(\lambda^{q-1}_{\max}\left(\mathcal{X}\#_f\mathcal{Y}\right)\mathcal{X}\#_f\mathcal{Y}\npreceq \mathcal{C} \right) 
\leq \mathrm{Tr}\left(\mathbb{E}\left[\left(\mathcal{X}^q \#_f \mathcal{Y}^q\right)^p\right] \star_N \mathcal{C}\right),
\end{eqnarray}
where $q $$\geq$$ 1$ and $p $$\geq$$ 1$, and 
\begin{eqnarray}\label{eq2-3-1:cor:3.3}
\mathrm{Pr}\left(\mathcal{X}^q \#_f \mathcal{Y}^q \npreceq \mathcal{C} \right) 
\leq \mathrm{Tr}\left(\mathbb{E}\left[\left(\lambda^{q-1}_{\max}\left(\mathcal{X}\#_f\mathcal{Y}\right)\mathcal{X}\#_f\mathcal{Y}\right)^p\right] \star_N \mathcal{C}\right),
\end{eqnarray}
where $0 $$\leq$$ q $$\leq$$ 1$ and $p $$\geq$$ 1$.

On the other hand, if the function $f$ is a pmd function, we have
\begin{eqnarray}\label{eq3-3:cor:3.3}
\mathrm{Pr}\left(\mathcal{X}^q \#_f \mathcal{Y}^q \npreceq \mathcal{C} \right) 
\leq \mathrm{Tr}\left(\mathbb{E}\left[\left( \lambda^{q-1}_{\min}\left(\mathcal{X}\#_f\mathcal{Y}\right)\mathcal{X}\#_f\mathcal{Y}\right)^p\right] \star_N \mathcal{C}\right),
\end{eqnarray}
where $q $$>$$ 1$ and $p $$\geq$$ 1$, and 
\begin{eqnarray}\label{eq3-3-1:cor:3.3}
\mathrm{Pr}\left(\lambda^{q-1}_{\min}\left(\mathcal{X}\#_f\mathcal{Y}\right)\mathcal{X}\#_f\mathcal{Y}\npreceq \mathcal{C} \right) 
\leq \mathrm{Tr}\left(\mathbb{E}\left[\left(\mathcal{X}^q \#_f \mathcal{Y}^q\right)^p\right] \star_N \mathcal{C}\right),
\end{eqnarray}
where $0 $$\leq$$ q $$\leq$$ 1$ and $p $$\geq$$ 1$.
\end{corollary} 
\textbf{Proof:}
If the function $f$ is a pmi function, we have $f(\mathcal{C}^q) $$\succeq$$ f^q(\mathcal{D})$ for any PD tensor $\mathcal{D}$. Then, we have 
\begin{eqnarray}\label{eq2:cor:3.3}
\Psi_{lower}\left(q,f,\mathcal{X},\mathcal{Y}\right) &\geq& 1,~\mbox{for $q $$\geq$$ 1$,} \nonumber \\
\lambda_{\min}\left(f^{-q}\left(\mathcal{Z}_0\right)f\left(\mathcal{Z}_0^q\right)\right)&\leq& 1,~\mbox{for $0 $$\leq$$ q $$\leq$$ 1$.}
\end{eqnarray}

For $q $$\geq$$ 1$, from Eq.~\eqref{eq2:cor:3.3} and Eq.~\eqref{eq5-3:thm:3.2}, we have 
\begin{eqnarray}\label{eq2-1:cor:3.3}
\lambda^{q-1}_{\max}\left(\mathcal{X}\#_f\mathcal{Y}\right)\mathcal{X}\#_f\mathcal{Y} &\preceq& \mathcal{X}^q \#_f \mathcal{Y}^q,
\end{eqnarray}
and for $0 $$< $$q $$\leq$$ 1$, from Eq.~\eqref{eq2:cor:3.3} and Eq.~\eqref{eq6-4:thm:3.2}, we have 
\begin{eqnarray}\label{eq2-2:cor:3.3}
\lambda^{q-1}_{\max}\left(\mathcal{X}\#_f\mathcal{Y}\right)\mathcal{X}\#_f\mathcal{Y} &\succeq& \mathcal{X}^q \#_f \mathcal{Y}^q.
\end{eqnarray}
Applying Lemma~\ref{lma:Loewner ordering with Markov Cheb inequalities} to Eq.~\eqref{eq2-1:cor:3.3} and Eq.~\eqref{eq2-2:cor:3.3}, we have the desired results at Eq.~\eqref{eq2-3:cor:3.3} and Eq.~\eqref{eq2-3-1:cor:3.3}.

If the function $f$ is a pmd function, we have $f(\mathcal{D}^q) $$\preceq$$ f^q(\mathcal{D})$ for any PD tensor $\mathcal{D}$. Then, we have 
\begin{eqnarray}\label{eq3:cor:3.3}
\Psi_{upper}\left(q,f,\mathcal{X},\mathcal{Y}\right) &\leq& 1,~\mbox{for $q $$\geq$$ 1$,} \nonumber \\
\lambda_{\max}\left(f^{-q}\left(\mathcal{Z}_0\right)f\left(\mathcal{Z}_0^q\right)\right)&\geq& 1,~\mbox{for $0 $$\leq$$ q $$\leq$$ 1$.}
\end{eqnarray}
For $q $$\geq$$ 1$, from Eq.~\eqref{eq3:cor:3.3} and Eq.~\eqref{eq5-3:thm:3.2}, we have 
\begin{eqnarray}\label{eq3-1:cor:3.3}
\mathcal{X}^q \#_f \mathcal{Y}^q &\preceq& \lambda^{q-1}_{\min}\left(\mathcal{X}\#_f\mathcal{Y}\right)\mathcal{X}\#_f\mathcal{Y},
\end{eqnarray}
and for $0 $$< $$q $$\leq$$ 1$, from Eq.~\eqref{eq3:cor:3.3} and Eq.~\eqref{eq6-4:thm:3.2}, we have 
\begin{eqnarray}\label{eq3-2:cor:3.3}
\mathcal{X}^q \#_f \mathcal{Y}^q &\succeq& \lambda^{q-1}_{\min}\left(\mathcal{X}\#_f\mathcal{Y}\right)\mathcal{X}\#_f\mathcal{Y}.
\end{eqnarray}
Applying Lemma~\ref{lma:Loewner ordering with Markov Cheb inequalities} tp Eq.~\eqref{eq3-1:cor:3.3} and Eq.~\eqref{eq3-2:cor:3.3}, we prove Eq.~\eqref{eq3-3:cor:3.3} and Eq.~\eqref{eq3-3-1:cor:3.3}.
$\hfill \Box$

Theorem~\ref{thm:3.2} is based on the function $f  $$\in$$ \mbox{TMI}^{1}$. Next theorem is to consider tail bounds for the function $h  $$\in$$ \mbox{TMD}^{1}$.
\begin{theorem}\label{thm:Prop. 3.9}
Given two random PD tensors $\mathcal{X} $$\in$$ \mathbb{C}^{I_1 \times \cdots \times I_N \times I_1 \times \cdots \times I_N}$, $\mathcal{Y} $$\in$$ \mathbb{C}^{I_1 \times \cdots \times I_N \times I_1 \times \cdots \times I_N}$ and a PD determinstic tensor $\mathcal{C}$, if $q=2^n q_0 \geq 1$ with $1 \leq q_0 \leq 2$, we set $\mathcal{Z}_{k-1} \define \mathcal{X}^{-2^{k-2}} \mathcal{Y}^{2^{k-1}}\mathcal{X}^{-2^{k-2}}$ for $k=1,2,\cdots,n$. We assume that $\mathcal{X} $$\#_h$$ \mathcal{Y} $$\preceq$$ \mathcal{I}$ almost surely with $h $$\in$$\mbox{TMD}^{1}$, we have 
\begin{eqnarray}\label{eq1-1:Prop. 3.9}
\mathrm{Pr}\left(\mathcal{X}^q \#_h \mathcal{Y}^q \npreceq \mathcal{C} \right)
&\leq& \mathrm{Tr}\left(\mathbb{E}\left[\left(\Phi_{upper}\left(q,h,\mathcal{X},\mathcal{Y}\right)\lambda^{q-1}_{\min}\left(\mathcal{X}\#_h\mathcal{Y}\right)\mathcal{X}\#_h\mathcal{Y}\right)^p\right]\star_N \mathcal{C}^{-1} \right),
\end{eqnarray}
and
\begin{eqnarray}\label{eq1-2:Prop. 3.9}
\mathrm{Pr}\left(\Phi_{lower}\left(q,h,\mathcal{X},\mathcal{Y}\right)\lambda^{q-1}_{\max}\left(\mathcal{X}\#_h\mathcal{Y}\right)\mathcal{X}\#_h\mathcal{Y} \npreceq \mathcal{C} \right)
&\leq& \mathrm{Tr}\left(\mathbb{E}\left[\left(\mathcal{X}^q \#_h \mathcal{Y}^q\right)^p\right]\star_N \mathcal{C}^{-1} \right),
\end{eqnarray}
where $\Phi_{lower}\left(q,h,\mathcal{X},\mathcal{Y}\right)$ and $\Phi_{upper}\left(q,h,\mathcal{X},\mathcal{Y}\right)$ are two positive numbers defined by
\begin{eqnarray}
\Phi_{lower}\left(q,h,\mathcal{X},\mathcal{Y}\right)&\define&\lambda_{\min}\left(h^{-q_0}\left(\mathcal{Z}_n\right)h\left(\mathcal{Z}^{q_0}_n\right)\right) \acute{\prod}_{k=1}^n \lambda_{\min}\left(h^{-2}\left(\mathcal{Z}_{k-1}\right)h\left(\mathcal{Z}_{k-1}^{2}\right)\right) \nonumber \\
\Phi_{upper}\left(q,h,\mathcal{X},\mathcal{Y}\right)&\define&\lambda_{\max}\left(h^{-q_0}\left(\mathcal{Z}_n\right)(h\left(\mathcal{Z}_n^{q_0}\right)\right)\acute{\prod}_{k=1}^n \lambda_{\max}\left(h^{-2}\left(\mathcal{Z}_{k-1}\right)h\left(\mathcal{Z}_{k-1}^{2}\right)\right).
\end{eqnarray}
Note that the definition of $\acute{\prod}$ is provided by Eq.~\eqref{eq4-1-2:thm:3.2}.

For $0 < q  \leq1$, we have 
\begin{eqnarray}\label{eq2-1:Prop. 3.9}
\mathrm{Pr}\left(\mathcal{X}^q \#_h \mathcal{Y}^q \npreceq \mathcal{C} \right)
&\leq& \mathrm{Tr}\left(\mathbb{E}\left[\left(\lambda_{\max}\left(h^{-q}\left(\mathcal{Z}_0\right)(h\left(\mathcal{Z}_0^q\right)\right)\lambda^{q-1}_{\min}\left(\mathcal{X}\#_h\mathcal{Y}\right)\mathcal{X}\#_h\mathcal{Y}\right)^p\right]\star_N \mathcal{C}^{-1} \right),
\end{eqnarray}
and
\begin{eqnarray}\label{eq2-2:Prop. 3.9}
\mathrm{Pr}\left( \lambda_{\min}\left(h^{-q}\left(\mathcal{Z}_0\right)h\left(\mathcal{Z}_0^q\right)\right)\lambda^{q-1}_{\max}\left(\mathcal{X}\#_h\mathcal{Y}\right)\mathcal{X}\#_h\mathcal{Y} \npreceq \mathcal{C} \right)
&\leq& \mathrm{Tr}\left(\mathbb{E}\left[\left(\mathcal{X}^q \#_h \mathcal{Y}^q\right)^p\right]\star_N \mathcal{C}^{-1} \right).
\end{eqnarray}
where $p \geq 1$. 
\end{theorem}
\textbf{Proof:}
As before, we will remove $\star_N$ in this proof for simplification. We begin with the case for $q \geq 1$.  We will separte the region of $q \geq 1$ into $1 \leq q \leq 2$ and $q \geq 2$. 

For the subregion $1 \leq q \leq 2$, $r \define 2-q$, and $\mathcal{X} $$\#_f$$ \mathcal{Y} $$\preceq$$ \mathcal{I}$, we have
\begin{eqnarray}\label{eq4:Prop. 3.9}
\mathcal{X}^q \#_h \mathcal{Y}^q &=& \mathcal{X}^{\frac{q}{2}}h\left(\mathcal{X}^{\frac{1-q}{2}}\mathcal{Z}_0\mathcal{X}^{\frac{1}{2}}\mathcal{Y}^{-r}\mathcal{X}^{\frac{1}{2}}\mathcal{Z}_0\mathcal{X}^{\frac{1-q}{2}}\right)\mathcal{X}^{\frac{q}{2}} \nonumber \\
&=& \mathcal{X}^{\frac{q}{2}}h\left(\mathcal{X}^{\frac{1-q}{2}}\mathcal{Z}_0\mathcal{X}^{\frac{1}{2}}\left(\mathcal{X}^{\frac{-1}{2}} \mathcal{Z}^{-1} \mathcal{X}^{\frac{-1}{2}}\right)^{r}\mathcal{X}^{\frac{1}{2}}\mathcal{Z}_0\mathcal{X}^{\frac{1-q}{2}}\right)\mathcal{X}^{\frac{q}{2}} \nonumber \\
&=& \mathcal{X}^{\frac{1}{2}}\mathcal{X}^{\frac{1-r}{2}}h\left(\mathcal{X}^{\frac{r-1}{2}}\mathcal{Z}_0
\left(\mathcal{X} \#_{x^r} \mathcal{Z}_)^{-1}\right)
\mathcal{Z}_0\mathcal{X}^{\frac{r-1}{2}}\right)\mathcal{X}^{\frac{1-r}{2}}\mathcal{X}^{\frac{1}{2}}\nonumber \\
&=& \mathcal{X}^{\frac{1}{2}}\left[\mathcal{X}^{1-r}\#_h \left(\mathcal{Z}_0\left(\mathcal{X} \#_{x^r} \mathcal{Z}_0^{-1}\right)\mathcal{Z}_0\right)\right]\mathcal{X}^{\frac{1}{2}}\nonumber \\
&\preceq_1&  \mathcal{X}^{\frac{1}{2}}\left[h^{1-q}(\mathcal{Z}_0)\#_h \left(\mathcal{Z}_0\left(h^{-1}(\mathcal{Z}_0)\#_{x^r} \mathcal{Z}_0^{-1}\right)\mathcal{Z}_0\right)\right]\mathcal{X}^{\frac{1}{2}}\nonumber \\
&=& \mathcal{X}^{\frac{1}{2}}\left[h^{1-q}\left(\mathcal{Z}_0\right)(h\left(\mathcal{Z}_0^q \right)\right]\mathcal{X}^{\frac{1}{2}}\nonumber \\
&\preceq& \lambda_{\max}\left(h^{-q}\left(\mathcal{Z}_0\right)(h\left(\mathcal{Z}_0^q\right)\right)\mathcal{X}\#_h\mathcal{Y},
\end{eqnarray}
where $\preceq_1$ we applies Lemma~\ref{lma:2.3} based on $\mathcal{X} $$\#_f$$ \mathcal{Y} $$\preceq$$ \mathcal{I} \Longleftrightarrow \mathcal{X} \preceq f^{-1}(\mathcal{Z})$. If $q \geq 2$, by finding some natural number $n$ such that $q=2^n q_0$ with $1 \leq q_0 \leq 2$ and iterating the relation provided by Eq.~\eqref{eq4:Prop. 3.9} dyadically with respect to $q$, we have the relation:
\begin{eqnarray}\label{eq4-1:Prop. 3.9}
\mathcal{X}^q \#_h \mathcal{Y}^q &\preceq& \lambda_{\max}\left(h^{-q_0}\left(\mathcal{Z}_n\right)h\left(\mathcal{Z}^{q_0}_n\right)\right) \prod_{k=1}^n \lambda_{\max}\left(h^{-2}\left(\mathcal{Z}_{k-1}\right)h\left(\mathcal{Z}_{k-1}^{2}\right)\right) \mathcal{X} \#_h \mathcal{Y}. 
\end{eqnarray}
By combining Eq.~\eqref{eq4:Prop. 3.9} and Eq.~\eqref{eq4-1:Prop. 3.9}, for $q \geq 1$, we have 
\begin{eqnarray}\label{eq4-1-1:Prop. 3.9}
\mathcal{X}^q \#_h \mathcal{Y}^q &\preceq& \lambda_{\max}\left(h^{-q_0}\left(\mathcal{Z}_n\right)h\left(\mathcal{Z}^{q_0}_n\right)\right) \acute{\prod}_{k=1}^n \lambda_{\max}\left(h^{-2}\left(\mathcal{Z}_{k-1}\right)h\left(\mathcal{Z}_{k-1}^{2}\right)\right) \mathcal{X} \#_h \mathcal{Y}.
\end{eqnarray}

Given any positive real number, say $\beta$, we can replace $\mathcal{X}$ and $\mathcal{Y}$ in Eq.~\eqref{eq4-1-1:Prop. 3.9} with 
$\beta^{-1}\mathcal{X}$ and $\beta^{-1}\mathcal{Y}$ to get
\begin{eqnarray}\label{eq4-2:Prop. 3.9}
\mathcal{X}^q \#_h \mathcal{Y}^q &\preceq& \lambda_{\max}\left(h^{-q_0}\left(\mathcal{Z}_n\right)h\left(\mathcal{Z}^{q_0}_n\right)\right) \acute{\prod}_{k=1}^n \lambda_{\max}\left(h^{-2}\left(\mathcal{Z}_{k-1}\right)h\left(\mathcal{Z}_{k-1}^{2}\right)\right)\beta^{q-1} \mathcal{X} \#_h \mathcal{Y}.
\end{eqnarray}
We can select $\beta $$=$$ \lambda_{\min}\left(\mathcal{X}\#_h\mathcal{Y}\right)$ to associate $\beta$ with 
$\mathcal{X}\#_h\mathcal{Y}$ and minimize R.H.S. of Eq.~\eqref{eq4-2:Prop. 3.9} in the sense of Loewner ordering. Then, we have 
\begin{eqnarray}\label{eq4-3:Prop. 3.9}
\mathcal{X}^q \#_f \mathcal{Y}^q\preceq\lambda_{\max}\left(f^{-q_0}\left(\mathcal{Z}_n\right)f\left(\mathcal{Z}^{q_0}_n\right)\right) \acute{\prod}_{k=1}^n \lambda_{\max}\left(h^{-2}\left(\mathcal{Z}_{k-1}\right)h\left(\mathcal{Z}_{k-1}^{2}\right)\right)\lambda^{q-1}_{\min}\left(\mathcal{X}\#_h\mathcal{Y}\right)\mathcal{X}\#_h\mathcal{Y}.
\end{eqnarray}

By replacing $h$, $\mathcal{X}$ and $\mathcal{Y}$ in Eq.~\eqref{eq4-3:Prop. 3.9} with  $h^{\ast} \define h^{-1}(x^{-1})$ for $x $$\in$$ (0,\infty)$, $\mathcal{X}^{-1}$ and $\mathcal{Y}^{-1}$, we have 
\begin{eqnarray}\label{eq5:Prop. 3.9}
\mathcal{X}^{-q} \#_{h^\ast} \mathcal{Y}^{-q} &\preceq& \lambda_{\max}\left((h^{\ast})^{-q_0}\left(\mathcal{Z}_n'\right)h^{\ast}\left(\mathcal{Z}_n'^{q_0}\right)\right)\acute{\prod}_{k=1}^n \lambda_{\max}\left((h^\ast)^{-2}\left(\mathcal{Z}'_{k-1}\right)(h^{\ast})\left(\mathcal{Z}_{k-1}^{'2}\right)\right)\nonumber \\
& & \cdot \lambda^{q-1}_{\min}\left(\mathcal{X}^{-1}\#_{h^{\ast}}\mathcal{Y}^{-1}\right)\mathcal{X}^{-1}\#_{h^{\ast}}\mathcal{Y}^{-1}.
\end{eqnarray}
where $\mathcal{Z}'_{k-1} \define \mathcal{X}^{2^{k-2}} \mathcal{Y}^{-2^{k-1}}\mathcal{X}^{2^{k-2}}$ for $k=1,2,\cdots,n$. From the definition of $h^{\ast}$, we have 
\begin{eqnarray}\label{eq5-1:Prop. 3.9}
\mathcal{X}^{-1}\#_{h^{\ast}}\mathcal{Y}^{-1}&=&\left(\mathcal{X}\#_{h}\mathcal{Y}\right)^{-1}, \nonumber \\
\mathcal{X}^{-q}\#_{h^{\ast}}\mathcal{Y}^{-q}&=&\left(\mathcal{X}^q \#_{h} \mathcal{Y}^q\right)^{-1},\nonumber \\
(h^{\ast})^{-q_0}\left(\mathcal{Z}_n'\right)h^{\ast}\left(\mathcal{Z}_n'^{q_0}\right)&=&h^{-1}\left(\mathcal{Z}_n^{q_0}\right)h^{q_0}\left(\mathcal{Z}_n\right),\nonumber \\
(h^{\ast})^{-2}\left(\mathcal{Z}_{k-1}'\right)h^{\ast}\left(\mathcal{Z}_{k-1}'^{2}\right)&=&h^{-1}\left(\mathcal{Z}_{k-1}^{2}\right)h^{2}\left(\mathcal{Z}_{k-1}\right),~\mbox{for $k=1,2,\cdots,n$;}
\end{eqnarray} 
then, by applying Eq.~\eqref{eq5-1:Prop. 3.9} to Eq.~\eqref{eq5:Prop. 3.9}, we obtain the following:
\begin{eqnarray}\label{eq5-2:Prop. 3.9}
\mathcal{X}^q \#_h \mathcal{Y}^q &\succeq& \lambda_{\min}\left(f^{-q_0}\left(\mathcal{Z}_n\right)(f\left(\mathcal{Z}_n^{q_0}\right)\right)\acute{\prod}_{k=1}^n \lambda_{\min}\left(f^{-2}\left(\mathcal{Z}_{k-1}\right)f\left(\mathcal{Z}_{k-1}^{2}\right)\right)\nonumber \\
&& \cdot \lambda^{q-1}_{\min}\left(\mathcal{X}\#_h\mathcal{Y}\right)\mathcal{X}\#_h\mathcal{Y}.
\end{eqnarray} 
By combining Eq.~\eqref{eq4-3:Prop. 3.9} and Eq.~\eqref{eq5-2:Prop. 3.9}, for $q \geq 1$, we have 
\begin{eqnarray}\label{eq5-3:Prop. 3.9}
\Phi_{lower}\left(q,h,\mathcal{X},\mathcal{Y}\right)\lambda^{q-1}_{\max}\left(\mathcal{X}\#_h\mathcal{Y}\right)\mathcal{X}\#_h\mathcal{Y} &\preceq& \mathcal{X}^q \#_h \mathcal{Y}^q  \nonumber \\
&\preceq &\Phi_{upper}\left(q,h,\mathcal{X},\mathcal{Y}\right)\lambda^{q-1}_{\min}\left(\mathcal{X}\#_h\mathcal{Y}\right)\mathcal{X}\#_h\mathcal{Y}.
\end{eqnarray}

By applying Lemma~\ref{lma:Loewner ordering with Markov Cheb inequalities}, we have the desired bounds provided by Eq.~\eqref{eq1-1:Prop. 3.9} and Eq.~\eqref{eq1-2:Prop. 3.9}.

Now, we will consider the case for $0 < q \leq 1$. 
\begin{eqnarray}\label{eq6:Prop. 3.9}
\mathcal{X}^q \#_h \mathcal{Y}^q
&=& \mathcal{X}^{\frac{q}{2}}h\left(\mathcal{X}^{\frac{-q}{2}}\left(\mathcal{X}^{\frac{1}{2}}\mathcal{Z}_0\mathcal{X}^{\frac{1}{2}}\right)^{q}\mathcal{X}^{\frac{-q}{2}}\right)\mathcal{X}^{\frac{q}{2}} \nonumber \\
&=& \mathcal{X}^{\frac{q}{2}}h\left(\mathcal{X}^{\frac{1-q}{2}}\left(\mathcal{X}^{-1}\#_{x^q}\mathcal{Z}_0\right)\mathcal{X}^{\frac{1-q}{2}}\right)\mathcal{X}^{\frac{q}{2}}\nonumber \\
&=& \mathcal{X}^{\frac{1}{2}}\left[\mathcal{X}^{-(1-q)}\#_h \left(\mathcal{X}^{-1}\#_{x^q}\mathcal{Z}_0\right)\right]\mathcal{X}^{\frac{1}{2}}\nonumber \\
&\succeq_1& \mathcal{X}^{\frac{1}{2}}\left[h^{(1-q)}\left(\mathcal{Z}_0\right)\#_h \left(h\left(\mathcal{Z}_0\right)\#_{x^q}\mathcal{Z}_0\right)\right]\mathcal{X}^{\frac{1}{2}}\nonumber \\
&=& \mathcal{X}^{\frac{1}{2}}\left[h^{(1-q)}\left(\mathcal{Z}_0\right)(h\left(\mathcal{Z}_0^q\right)\right]\mathcal{X}^{\frac{1}{2}}\nonumber \\
&\succeq& \lambda_{\min}\left(h^{-q}\left(\mathcal{Z}_0\right)(h\left(\mathcal{Z}_0^q\right)\right)\mathcal{X}\#_h\mathcal{Y},
\end{eqnarray}
where $\succeq_1$ we applies  Lemma~\ref{lma:2.3} based on $\mathcal{X} $$\#_h$$ \mathcal{Y} $$\preceq$$ \mathcal{I} \Longleftrightarrow \mathcal{X} \preceq h^{-1}(\mathcal{Z}_0)$ with $0 \leq 1-q < 1$. Given any positive real number, say $\beta$, we also can replace $\mathcal{X}$ and $\mathcal{Y}$ in Eq.~\eqref{eq6:Prop. 3.9} with 
$\beta^{-1}\mathcal{X}$ and $\beta^{-1}\mathcal{Y}$ to get
\begin{eqnarray}\label{eq6-1:Prop. 3.9}
\mathcal{X}^q \#_f \mathcal{Y}^q &\succeq& \lambda_{\min}\left(h^{-q}\left(\mathcal{Z}_0\right)(h\left(\mathcal{Z}_0^q\right)\right)\beta^{q-1}\mathcal{X}\#_f\mathcal{Y}.
\end{eqnarray}
We can select $\beta $$=$$ \lambda_{\max}\left(\mathcal{X}\#_h\mathcal{Y}\right)$ to associate $\beta$ with 
$\mathcal{X}\#_h\mathcal{Y}$ and maximize R.H.S. of Eq.~\eqref{eq6-1:Prop. 3.9} in the sense of L\"oewner ordering. Then, we have 
\begin{eqnarray}\label{eq6-2:Prop. 3.9}
\mathcal{X}^q \#_h \mathcal{Y}^q &\succeq& \lambda_{\min}\left(h^{-q}\left(\mathcal{Z}_0\right)h\left(\mathcal{Z}_0^q\right)\right)\lambda^{q-1}_{\max}\left(\mathcal{X}\#_h\mathcal{Y}\right)\mathcal{X}\#_h\mathcal{Y}.
\end{eqnarray}

By replacing $h $$\in$$ \mbox{TMD}^{1}$, $\mathcal{X}$ and $\mathcal{Y}$ in Eq.~\eqref{eq6-2:Prop. 3.9} with  $h^{\ast} \define h^{-1}(x^{-1})$ for $x $$\in$$ (0,\infty)$, $\mathcal{X}^{-1}$ and $\mathcal{Y}^{-1}$, we have 
\begin{eqnarray}\label{eq6-3:Prop. 3.9}
\mathcal{X}^q \#_h \mathcal{Y}^q &\preceq& \lambda_{\max}\left(h^{-q}\left(\mathcal{Z}_0\right)h\left(\mathcal{Z}_0^q\right)\right)\lambda^{q-1}_{\min}\left(\mathcal{X}\#_h\mathcal{Y}\right)\mathcal{X}\#_h \mathcal{Y}.
\end{eqnarray}
Therefore, combining Eq.~\eqref{eq6-2:Prop. 3.9} and Eq.~\eqref{eq6-3:Prop. 3.9}, we have 
\begin{eqnarray}\label{eq6-4:Prop. 3.9}
\lambda_{\min}\left(h^{-q}\left(\mathcal{Z}_0\right)h\left(\mathcal{Z}_0^q\right)\right)\lambda^{q-1}_{\max}\left(\mathcal{X}\#_h\mathcal{Y}\right)\mathcal{X}\#_h\mathcal{Y} &\preceq& \mathcal{X}^q \#_h \mathcal{Y}^q  \nonumber \\
&\preceq &  \lambda_{\max}\left(h^{-q}\left(\mathcal{Z}_0\right)(h\left(\mathcal{Z}_0^q\right)\right)\lambda^{q-1}_{\min}\left(\mathcal{X}\#_h\mathcal{Y}\right)\mathcal{X}\#_h\mathcal{Y},\nonumber \\
\end{eqnarray}
where $0 \leq q \leq 1$. By applying Lemma~\ref{lma:Loewner ordering with Markov Cheb inequalities}, we have the desired bounds provided by Eq.~\eqref{eq2-1:Prop. 3.9} and Eq.~\eqref{eq2-2:Prop. 3.9}.
$\hfill \Box$

\subsection{Connection Functions Come From $\mbox{TC}^{1}$}\label{sec:Connection Functions are TC1}

In this section, we will consider connection functions come from $\mbox{TC}^{1}$. Before presenting those main results in this section, we will present the following Lemma about Kantorovich type inequality for operators~\cite{furuta2005mond}.
\begin{lemma}\label{lma:Kantorovich type inequality}
Let $\mathcal{A}, \mathcal{B} $$\in$$ \mathbb{C}^{I_1 \times \cdots \times I_N \times I_1 \times \cdots \times I_N}$ be two PD tensors such that 
\begin{eqnarray}\label{eq1:lma:Kantorovich type inequality}
m_1 \mathcal{I} \preceq \mathcal{A} \preceq M_1 \mathcal{I},~\mbox{and} \nonumber \\
m_2 \mathcal{I} \preceq \mathcal{B} \preceq M_2 \mathcal{I},~~~~~~~
\end{eqnarray}
where $M_1 $$>$$ m_1 $$>$$ 0$, and $M_2 $$>$$ m_2 $$>$$ 0$. If $\mathcal{B} $$\preceq$$ \mathcal{A}$ and $p $$>$$ 1$, we have
\begin{eqnarray}\label{eq2:lma:Kantorovich type inequality}
\mathcal{B}^p &\preceq& \mathrm{K}(m_1,M_1,p) \mathcal{A}^p, \nonumber \\
\mathcal{B}^p &\preceq& \mathrm{K}(m_2,M_2,p) \mathcal{A}^p,
\end{eqnarray}
where the Kantorovich contant, $\mathrm{K}(m,M,p)$, can be expressed by
\begin{eqnarray}
\mathrm{K}(m,M,p) &=& \left(\frac{(p-1)\left(M^p - m^p\right)}{p\left(mM^p - Mm^p\right)}\right)^p\frac{mM^p - Mm^p}{(p-1)(M-m)}. 
\end{eqnarray}
\end{lemma}
\textbf{Proof:}
Theorem 8.3 from~\cite{furuta2005mond}.
$\hfill \Box$

Following theorem is about the tail bounds for connection functions come from $\mbox{TC}^{1}$.
\begin{theorem}\label{thm:Prop. 3.10}
Given two random PD tensors $\mathcal{X} $$\in$$ \mathbb{C}^{I_1 \times \cdots \times I_N \times I_1 \times \cdots \times I_N}$, $\mathcal{Y} $$\in$$ \mathbb{C}^{I_1 \times \cdots \times I_N \times I_1 \times \cdots \times I_N}$, and a PD determinstic tensor $\mathcal{C} $$\in$$ \mathbb{C}^{I_1 \times \cdots \times I_N \times I_1 \times \cdots \times I_N}$, we will set $\mathcal{Z} \define \mathcal{X}^{1/2}\star_N\mathcal{Y}^{-1}\star_N\mathcal{X}^{1/2}$. Let $g $$\in$$\mbox{TC}^{1}$, if $\mathcal{X} \#_g \mathcal{Y} \preceq \mathcal{I}$ almost surely, and $p,q \geq 1$, we have 
\begin{eqnarray}\label{eq1-1:Prop. 3.10}
\mathrm{Pr}\left(\mathcal{X}^q \#_g \mathcal{Y}^q \npreceq \mathcal{C}\right) \leq
\mathrm{Tr}\left(\mathbb{E}\left[\left(\mathrm{K}_1 \lambda^{1-q}_{\min}\left(\mathcal{X}\#_{g}\mathcal{Y}\right)
\lambda_{\max}\left(g^{-q}(\mathcal{Z})g(\mathcal{Z}^q)\right)\mathrm{K}_2\mathcal{I}\right)^p\right]\star_N \mathcal{C}^{-1}\right)
\end{eqnarray}
where $\mathrm{K}_1$ and $\mathrm{K}_2$ are set as
\begin{eqnarray}
\mathrm{K}_1 &\define& \mathrm{K}\left(\lambda^{-1}_{\max}\left(\mathcal{X}\right),\lambda^{-1}_{\min}\left(\mathcal{X}\right),q-1\right) \nonumber \\
\mathrm{K}_2 &\define& \mathrm{K}\left(\lambda^{-1}_{\max}\left(\mathcal{X}\right),\lambda^{-1}_{\min}\left(\mathcal{X}\right),2q-1\right).
\end{eqnarray}
Moreover, if $\mathcal{X} \#_g \mathcal{Y} \succeq \mathcal{I}$ almost surely, we have 
\begin{eqnarray}\label{eq1-2:Prop. 3.10}
\mathrm{Pr}\left(\lambda^{1-q}_{\min}\left(\mathcal{X}\#_{g}\mathcal{Y}\right)
\lambda_{\max}\left(g^{-q}(\mathcal{Z})g(\mathcal{Z}^q)\right) \mathrm{K}_2^{-1}\mathcal{I}\npreceq \mathcal{C}\right) \leq \mathrm{Tr}\left(\mathbb{E}\left[\left(\mathcal{X}^q \#_g \mathcal{Y}^q\right)^p\right] \star_N \mathcal{C}^{-1}\right)
\end{eqnarray}
\end{theorem}
\textbf{Proof:}
In this proof, we will remove $\star_N$ for presentation simplification. If we define $f(x)$ as $f(x) \define g(x)/x$ for $x >0$, we have 
\begin{eqnarray}
\mathcal{X}\#_g\mathcal{Y}&=&\mathcal{X}^{1/2}f(\mathcal{Z})\mathcal{X}^{1/2}.
\end{eqnarray}
We will prove Eq.~\eqref{eq1-1:Prop. 3.10} first. Given $\mathcal{X} \#_g \mathcal{Y} \preceq \mathcal{I}$, we have $f(\mathcal{Z}) $$\preceq$$ \mathcal{X}^{-1}$, almost surely. Since we have 
\begin{eqnarray}\label{eq2:Prop. 3.10}
\mathcal{X}^q \#_g \mathcal{Y}^q &=& \mathcal{X}^{q/2}f\left(\mathcal{X}^{q/2}\left(\mathcal{X}^{1/2}\mathcal{Z}^{-1}\mathcal{X}^{1/2}\right)^{-q}\mathcal{X}^{q/2}\right)\mathcal{X}^{q/2} \nonumber \\
&=& \mathcal{X}^{\frac{q}{2}}f\left(\mathcal{X}^{\frac{q-1}{2}}\mathcal{Z}\mathcal{X}^{-\frac{1}{2}}\left(\mathcal{X}^{\frac{1}{2}} \mathcal{Z}^{-1} \mathcal{X}^{\frac{1}{2}}\right)^{2-q}\mathcal{X}^{-\frac{1}{2}}\mathcal{Z}\mathcal{X}^{\frac{q-1}{2}}\right)\mathcal{X}^{\frac{q}{2}} \nonumber \\
&=& \mathcal{X}^{\frac{q}{2}}f\left(\mathcal{X}^{\frac{q-1}{2}}\mathcal{Z}
\left(\mathcal{X}^{-1} \#_{x^{2-q}} \mathcal{Z}^{-1}\right)
\mathcal{Z}\mathcal{X}^{\frac{q-1}{2}}\right)\mathcal{X}^{\frac{q}{2}}\nonumber \\
&=& \mathcal{X}^{q-\frac{1}{2}}\left(\mathcal{X}^{1-q}\#_f \left[\left(\mathcal{Z}\mathcal{X}^{-1}\mathcal{Z}\right)\#_{x^{2-q}}\mathcal{Z}\right]\right)\mathcal{X}^{q-\frac{1}{2}},
\end{eqnarray}
we will bound Eq.~\eqref{eq2:Prop. 3.10} for different value range of $q$.

Since
\begin{eqnarray}\label{eq3-1:Prop. 3.10}
\mathcal{X}^{-1} &\preceq& \lambda_{\max}\left(f^{-1/2}(\mathcal{Z})\mathcal{X}^{-1}f^{-1/2}(\mathcal{Z})\right)f(\mathcal{Z}) \nonumber \\
 &=& \lambda^{-1}_{\min}\left(\mathcal{X}\#_{g}\mathcal{Y}\right)f(\mathcal{Z}),
\end{eqnarray}
from Lemma~\ref{lma:2.3}, we will have
\begin{eqnarray}\label{eq3-2:Prop. 3.10}
\mathcal{X}^{1-q} &\preceq& \left(\lambda^{-1}_{\min}\left(\mathcal{X}\#_{g}\mathcal{Y}\right)f(\mathcal{Z})\right)^{q-1},
\end{eqnarray}
where $1 \leq $$q$$ \leq 2$.  

Note that $\lambda^{-1}_{\max}\left(\mathcal{X}\right)\mathcal{I} $$\preceq$$ \mathcal{X}^{-1} $$\preceq$$ \lambda^{-1}_{\min}\left(\mathcal{X}\right)\mathcal{I}$. Given $q \geq 2$ and Lemma~\ref{lma:Kantorovich type inequality}, Eq.~\eqref{eq3-1:Prop. 3.10} can be extended as
\begin{eqnarray}\label{eq3-3:Prop. 3.10}
\mathcal{X}^{1-q} &\preceq& \mathrm{K}\left(\lambda^{-1}_{\max}\left(\mathcal{X}\right),\lambda^{-1}_{\min}\left(\mathcal{X}\right),q-1\right)\left(\lambda^{-1}_{\min}\left(\mathcal{X}\#_{g}\mathcal{Y}\right)f(\mathcal{Z})\right)^{q-1}.
\end{eqnarray}
Because the Kantorovich contant $\mathrm{K}\left(\lambda^{-1}_{\max}\left(\mathcal{X}\right),\lambda^{-1}_{\min}\left(\mathcal{X}\right),q-1\right)$ is greater than $1$, we have
\begin{eqnarray}\label{eq3-4:Prop. 3.10}
\mathcal{X}^{1-q} &\preceq& \mathrm{K}\left(\lambda^{-1}_{\max}\left(\mathcal{X}\right),\lambda^{-1}_{\min}\left(\mathcal{X}\right),q-1\right)\left(\lambda^{-1}_{\min}\left(\mathcal{X}\#_{g}\mathcal{Y}\right)f(\mathcal{Z})\right)^{q-1},
\end{eqnarray}
for all $q \geq 1$. Also, from Eq.~\eqref{eq3-1:Prop. 3.10}, we have 
\begin{eqnarray}\label{eq3-5:Prop. 3.10}
\mathcal{Z}\mathcal{X}^{-1}\mathcal{Z} &\preceq& \lambda^{-1}_{\min}\left(\mathcal{X}\#_{g}\mathcal{Y}\right)\mathcal{Z}^2f(\mathcal{Z}).
\end{eqnarray}

At this status, we can upper bound Eq.~\eqref{eq2:Prop. 3.10} via L\"owner ordering as
\begin{eqnarray}\label{eq4-1:Prop. 3.10}
\mathcal{X}^q \#_g \mathcal{Y}^q &\preceq_1& \mathcal{X}^{q-\frac{1}{2}}\left(\mathrm{K}\left(\lambda^{-1}_{\max}\left(\mathcal{X}\right),\lambda^{-1}_{\min}\left(\mathcal{X}\right),q-1\right)\left(\lambda^{-1}_{\min}\left(\mathcal{X}\#_{g}\mathcal{Y}\right)f(\mathcal{Z})\right)^{q-1}\#_f \right. \nonumber \\
&&\left.  \left[\left(\lambda^{-1}_{\min}\left(\mathcal{X}\#_{g}\mathcal{Y}\right)\mathcal{Z}^2f(\mathcal{Z})\right)\#_{x^{2-q}}\mathcal{Z}\right]\right)\mathcal{X}^{q-\frac{1}{2}} \nonumber \\
&=& \mathrm{K}\left(\lambda^{-1}_{\max}\left(\mathcal{X}\right),\lambda^{-1}_{\min}\left(\mathcal{X}\right),q-1\right)\lambda^{1-q}_{\min}\left(\mathcal{X}\#_{g}\mathcal{Y}\right)
\mathcal{X}^{q-\frac{1}{2}}f^{q-1}(\mathcal{Z})f(\mathcal{Z}^q)\mathcal{X}^{q-\frac{1}{2}}  \nonumber \\
&=& \mathrm{K}\left(\lambda^{-1}_{\max}\left(\mathcal{X}\right),\lambda^{-1}_{\min}\left(\mathcal{X}\right),q-1\right)\lambda^{1-q}_{\min}\left(\mathcal{X}\#_{g}\mathcal{Y}\right)
\mathcal{X}^{q-\frac{1}{2}}\nonumber \\
&  &\star_N \left(g^{-q}(\mathcal{Z})g(\mathcal{Z}^q)f^{2q-1}(\mathcal{Z})\right)\mathcal{X}^{q-\frac{1}{2}}  \nonumber \\
&\preceq& \mathrm{K}\left(\lambda^{-1}_{\max}\left(\mathcal{X}\right),\lambda^{-1}_{\min}\left(\mathcal{X}\right),q-1\right)\lambda^{1-q}_{\min}\left(\mathcal{X}\#_{g}\mathcal{Y}\right)
\lambda_{\max}\left(g^{-q}(\mathcal{Z})g(\mathcal{Z}^q)\right)\nonumber \\
&  &\mathcal{X}^{q-\frac{1}{2}}f^{2q-1}(\mathcal{Z})\mathcal{X}^{q-\frac{1}{2}}\nonumber \\
&\preceq_2& \mathrm{K}\left(\lambda^{-1}_{\max}\left(\mathcal{X}\right),\lambda^{-1}_{\min}\left(\mathcal{X}\right),q-1\right)\lambda^{1-q}_{\min}\left(\mathcal{X}\#_{g}\mathcal{Y}\right)
\lambda_{\max}\left(g^{-q}(\mathcal{Z})g(\mathcal{Z}^q)\right) \nonumber \\
&  &\mathrm{K}\left(\lambda^{-1}_{\max}\left(\mathcal{X}\right),\lambda^{-1}_{\min}\left(\mathcal{X}\right),2q-1\right)\mathcal{I},
\end{eqnarray}
where we apply Eq.~\eqref{eq3-4:Prop. 3.10} and Eq.~\eqref{eq3-5:Prop. 3.10} at $\preceq_1$, and we apply Lemma~\ref{lma:Kantorovich type inequality} to $f(\mathcal{Z}) $$\preceq$$ \mathcal{X}^{-1}$ at $\preceq_2$. 
Eq.~\eqref{eq1-1:Prop. 3.10} is obtained from applying Lemma~\ref{lma:Loewner ordering with Markov Cheb inequalities} to Eq.~\eqref{eq4-1:Prop. 3.10}.

Now, we will prove Eq.~\eqref{eq1-2:Prop. 3.10}. Because, if $\mathcal{X} \#_g \mathcal{Y} \succeq \mathcal{I}$ almost surely and $1 \leq q \leq 2$,  we will have
\begin{eqnarray}\label{eq6-2:Prop. 3.10}
\mathcal{X}^{1-q} &\succeq& \left(\lambda^{-1}_{\max}\left(\mathcal{X}\#_{g}\mathcal{Y}\right)f(\mathcal{Z})\right)^{q-1}.
\end{eqnarray}
Given $q \geq 2$, due to $\lambda^{-1}_{\max}\left(\mathcal{X}\right)\mathcal{I} $$\preceq$$ \mathcal{X}^{-1} $$\preceq$$ \lambda^{-1}_{\min}\left(\mathcal{X}\right)\mathcal{I}$ and Lemma~\ref{lma:Kantorovich type inequality}, Eq.~\eqref{eq6-2:Prop. 3.10} can be extended as
\begin{eqnarray}\label{eq6-3:Prop. 3.10}
\mathrm{K}\left(\lambda^{-1}_{\max}\left(\mathcal{X}\right),\lambda^{-1}_{\min}\left(\mathcal{X}\right),q-1\right) \mathcal{X}^{1-q} &\succeq& \left(\lambda^{-1}_{\max}\left(\mathcal{X}\#_{g}\mathcal{Y}\right)f(\mathcal{Z})\right)^{q-1}.
\end{eqnarray}
Because the Kantorovich contant $\mathrm{K}\left(\lambda^{-1}_{\max}\left(\mathcal{X}\right),\lambda^{-1}_{\min}\left(\mathcal{X}\right),q-1\right)$ is greater than $1$, we have
\begin{eqnarray}\label{eq6-4:Prop. 3.10}
\mathcal{X}^{1-q} &\succeq& \left(\lambda^{-1}_{\max}\left(\mathcal{X}\#_{g}\mathcal{Y}\right)f(\mathcal{Z})\right)^{q-1},
\end{eqnarray}
for all $q \geq 1$. Also, from $\mathcal{X}^{-1} $$\succeq$$ \lambda^{-1}_{\max}\left(\mathcal{X}\#_{f}\mathcal{Y}\right)f(\mathcal{Z})$, we have 
\begin{eqnarray}\label{eq6-5:Prop. 3.10}
\mathcal{Z}\mathcal{X}^{-1}\mathcal{Z} &\succeq& \lambda^{-1}_{\max}\left(\mathcal{X}\#_{g}\mathcal{Y}\right)\mathcal{Z}^2f(\mathcal{Z}).
\end{eqnarray}

From Eq.~\eqref{eq2:Prop. 3.10}, we have
\begin{eqnarray}\label{eq7-1:Prop. 3.10}
\mathcal{X}^q \#_g \mathcal{Y}^q &\succeq_1& \mathcal{X}^{q-\frac{1}{2}}\left(\left(\lambda^{-1}_{\min}\left(\mathcal{X}\#_{g}\mathcal{Y}\right)f(\mathcal{Z})\right)^{q-1}\#_f  \left[\left(\lambda^{-1}_{\min}\left(\mathcal{X}\#_{g}\mathcal{Y}\right)\mathcal{Z}^2f(\mathcal{Z})\right)\#_{x^{2-q}}\mathcal{Z}\right]\right)\mathcal{X}^{q-\frac{1}{2}} \nonumber \\
&=& \lambda^{1-q}_{\min}\left(\mathcal{X}\#_{g}\mathcal{Y}\right)
\mathcal{X}^{q-\frac{1}{2}}f^{q-1}(\mathcal{Z})f(\mathcal{Z}^q)\mathcal{X}^{q-\frac{1}{2}}  \nonumber \\
&=& \lambda^{1-q}_{\min}\left(\mathcal{X}\#_{g}\mathcal{Y}\right)
\mathcal{X}^{q-\frac{1}{2}}\left(g^{-q}(\mathcal{Z})g(\mathcal{Z}^q)f^{2q-1}(\mathcal{Z})\right)\mathcal{X}^{q-\frac{1}{2}}  \nonumber \\
&\succeq& \lambda^{1-q}_{\min}\left(\mathcal{X}\#_{g}\mathcal{Y}\right)
\lambda_{\min}\left(g^{-q}(\mathcal{Z})g(\mathcal{Z}^q)\right)\mathcal{X}^{q-\frac{1}{2}}f^{2q-1}(\mathcal{Z})\mathcal{X}^{q-\frac{1}{2}}\nonumber \\
&\succeq_2& \lambda^{1-q}_{\min}\left(\mathcal{X}\#_{g}\mathcal{Y}\right)
\lambda_{\max}\left(g^{-q}(\mathcal{Z})g(\mathcal{Z}^q)\right) \mathrm{K}^{-1}\left(\lambda^{-1}_{\max}\left(\mathcal{X}\right),\lambda^{-1}_{\min}\left(\mathcal{X}\right),2q-1\right)\mathcal{I},
\end{eqnarray}
where we apply Eq.~\eqref{eq6-4:Prop. 3.10} and Eq.~\eqref{eq6-5:Prop. 3.10} at $\succeq_1$, and we apply Lemma~\ref{lma:Kantorovich type inequality} to $f(\mathcal{Z}) $$\succeq$$ \mathcal{X}^{-1}$ at $\succeq_2$. 

By applying Lemma~\ref{lma:Loewner ordering with Markov Cheb inequalities} to Eq.~\eqref{eq7-1:Prop. 3.10}, we can have the tail bound given by Eq.~\eqref{eq1-2:Prop. 3.10} for any $p $$\geq$$ 1$.
$\hfill \Box$

Following Corollary is based on Theorem~\ref{thm:Prop. 3.10} by transforming the function of $g$.
\begin{corollary}\label{cor:Prop. 3.10}
Given two random PD tensors $\mathcal{X} $$\in$$ \mathbb{C}^{I_1 \times \cdots \times I_N \times I_1 \times \cdots \times I_N}$, $\mathcal{Y} $$\in$$ \mathbb{C}^{I_1 \times \cdots \times I_N \times I_1 \times \cdots \times I_N}$, and a PD determinstic tensor $\mathcal{C} $$\in$$ \mathbb{C}^{I_1 \times \cdots \times I_N \times I_1 \times \cdots \times I_N}$, we will set $\mathcal{Z} \define \mathcal{X}^{1/2}\star_N\mathcal{Y}^{-1}\star_N\mathcal{X}^{1/2}$. Let $h $$\in$$\mbox{TMD}^{1}$, if $\mathcal{X} \#_h \mathcal{Y} \preceq \mathcal{I}$ almost surely, and $p,q \geq 1$, we have 
\begin{eqnarray}\label{eq1-1:cor:Prop. 3.10}
\mathrm{Pr}\left(\mathcal{X}^q \#_h \mathcal{Y}^q \npreceq \mathcal{C}\right) \leq
\mathrm{Tr}\left(\mathbb{E}\left[\left(\mathrm{K}_1 \lambda^{1-q}_{\min}\left(\mathcal{X}\#_{h}\mathcal{Y}\right)
\lambda_{\max}\left(h^{-q}(\mathcal{Z})h(\mathcal{Z}^q)\right)\mathrm{K}_2\mathcal{I}\right)^p\right]\star_N \mathcal{C}^{-1}\right)
\end{eqnarray}
where $\mathrm{K}_1$ and $\mathrm{K}_2$ are set as
\begin{eqnarray}
\mathrm{K}_1 &\define& \mathrm{K}\left(\lambda^{-1}_{\max}\left(\mathcal{Y}\right),\lambda^{-1}_{\min}\left(\mathcal{Y}\right),q-1\right) \nonumber \\
\mathrm{K}_2 &\define& \mathrm{K}\left(\lambda^{-1}_{\max}\left(\mathcal{Y}\right),\lambda^{-1}_{\min}\left(\mathcal{Y}\right),2q-1\right).
\end{eqnarray}
Moreover, if $\mathcal{X} \#_h \mathcal{Y} \succeq \mathcal{I}$ almost surely, we have 
\begin{eqnarray}\label{eq1-2:cor:Prop. 3.10}
\mathrm{Pr}\left(\lambda^{1-q}_{\min}\left(\mathcal{X}\#_{h}\mathcal{Y}\right)
\lambda_{\max}\left(h^{-q}(\mathcal{Z})h(\mathcal{Z}^q)\right) \mathrm{K}_2^{-1}\mathcal{I}\npreceq \mathcal{C}\right) \leq \mathrm{Tr}\left(\mathbb{E}\left[\left(\mathcal{X}^q \#_h \mathcal{Y}^q\right)^p\right] \star_N \mathcal{C}^{-1}\right).
\end{eqnarray}
\end{corollary}
\textbf{Proof:}
If the function $h(x)$ is expressed by $h(x) \define xg(x^{-1})$, then we have 
\begin{eqnarray}\label{eq5-0:cor:Prop. 3.10}
\mathcal{X}\#_g \mathcal{Y} = \mathcal{Y}\#_h \mathcal{X}.
\end{eqnarray}
\begin{eqnarray}\label{eq5-1:cor:Prop. 3.10}
g^{-q}(\mathcal{Z})g(\mathcal{Z}^q) = (\mathcal{Z}g(\mathcal{Z}^{-1}))^{-q}\mathcal{Z}^{q}g(\mathcal{Z}^{-q}) = h^{-q}(\mathcal{Z}^{-1})g(\mathcal{Z}^{-q}).
\end{eqnarray}
Moreover, we also have
\begin{eqnarray}\label{eq5-2:cor:Prop. 3.10}
\lambda_{\max}\left(g^{-q}(\mathcal{Z})g(\mathcal{Z}^q)\right) &=& \lambda_{\max}\left(h^{-q}(\mathcal{Y}^{1/2}\mathcal{X}^{-1}\mathcal{Y}^{1/2})g^q(\mathcal{Y}^{1/2}\mathcal{X}^{-1}\mathcal{Y}^{1/2})\right)\nonumber \\
&=& \lambda_{\max}\left( h^{-q}(\mathcal{X}^{-1/2}\mathcal{Y}\mathcal{X}^{-1/2})g^q(\mathcal{X}^{-1/2}\mathcal{Y}\mathcal{X}^{-1/2})\right);
\end{eqnarray}
and
\begin{eqnarray}\label{eq5-3:cor:Prop. 3.10}
\lambda_{\min}\left(g^{-q}(\mathcal{Z})g(\mathcal{Z}^q)\right) &=& \lambda_{\min}\left(h^{-q}(\mathcal{Y}^{1/2}\mathcal{X}^{-1}\mathcal{Y}^{1/2})g^q(\mathcal{Y}^{1/2}\mathcal{X}^{-1}\mathcal{Y}^{1/2})\right)\nonumber \\
&=& \lambda_{\min}\left( h^{-q}(\mathcal{X}^{-1/2}\mathcal{Y}\mathcal{X}^{-1/2})g^q(\mathcal{X}^{-1/2}\mathcal{Y}\mathcal{X}^{-1/2})\right).
\end{eqnarray}
Therefore, this Corollary is proved from Theorem~\ref{thm:Prop. 3.10} by using Eqs.~\eqref{eq5-0:cor:Prop. 3.10},~\eqref{eq5-1:cor:Prop. 3.10},~\eqref{eq5-2:cor:Prop. 3.10},~\eqref{eq5-3:cor:Prop. 3.10}.
$\hfill \Box$

\subsection{Bounds for $\Psi_{upper}$ and $\Psi_{lower}$ certain $f$}\label{sec:Bounds for Psi upper and lower}

Recall that $\Psi_{lower}\left(q,f,\mathcal{X},\mathcal{Y}\right)$ and $\Psi_{upper}\left(q,f,\mathcal{X},\mathcal{Y}\right)$ are two positive numbers defined by
\begin{eqnarray}\label{eq:Bounds for Psi upper lower}
\Psi_{lower}\left(q,f,\mathcal{X},\mathcal{Y}\right)&\define&\lambda_{\min}\left(f^{-q_0}\left(\mathcal{Z}_n\right)f\left(\mathcal{Z}^{q_0}_n\right)\right) \acute{\prod}_{k=1}^n \lambda_{\min}\left(f^{-2}\left(\mathcal{Z}_{k-1}\right)f\left(\mathcal{Z}_{k-1}^{2}\right)\right) \nonumber \\
\Psi_{upper}\left(q,f,\mathcal{X},\mathcal{Y}\right)&\define&\lambda_{\max}\left(f^{-q_0}\left(\mathcal{Z}_n\right)(f\left(\mathcal{Z}_n^{q_0}\right)\right)\acute{\prod}_{k=1}^n \lambda_{\max}\left(f^{-2}\left(\mathcal{Z}_{k-1}\right)f\left(\mathcal{Z}_{k-1}^{2}\right)\right).
\end{eqnarray}
where $q=2^n q_0 $$\geq$$ 1$. 

If the function $f$ satisfies that $\log f(e^x)$ is a convex function on $x \in (-\infty, \infty)$, we have the following bounds estimation Lemma for $\Psi_{upper}$ and $\Psi_{lower}$. We also define the following term for later notation simplicity:
\begin{eqnarray}
\psi(q,f,\mathcal{Z})&\define&\max\left(\frac{f(\lambda_{\min}^q(\mathcal{Z}))}{f^q(\lambda_{\min}(\mathcal{Z}))},\frac{f(\lambda_{\max}^q(\mathcal{Z}))}{f^q(\lambda_{\max}(\mathcal{Z}))}\right).
\end{eqnarray}

\begin{lemma}\label{lma:Bounds for Psi}
Given the function $f$ satisfing that $\log f(e^x)$ is a convex function on $x \in (-\infty, \infty)$, we have 
\begin{eqnarray}
\Psi_{lower}\left(q,f,\mathcal{X},\mathcal{Y}\right)&\geq&1.
\end{eqnarray}
For $\Psi_{upper}\left(q,f,\mathcal{X},\mathcal{Y}\right)$, we have 
\begin{eqnarray}
\Psi_{upper}\left(q,f,\mathcal{X},\mathcal{Y}\right)&\leq&\psi(q_0,f,\mathcal{Z}_n)
\acute{\prod}_{k=1}^n \psi(2,f,\mathcal{Z}_{k-1}),
\end{eqnarray}
where $\mathcal{Z}_{k-1} \define \mathcal{X}^{-2^{k-2}} \mathcal{Y}^{2^{k-1}}\mathcal{X}^{-2^{k-2}}$. 
\end{lemma}
\textbf{Proof:}
Since the function $\log f(e^x)$ is a convex function on $x \in (-\infty, \infty)$, we have 
\begin{eqnarray}\label{eq1:lma:Bounds for Psi}
\frac{\log f(e^x)}{dx} \geq 0.
\end{eqnarray}
Then, we have 
\begin{eqnarray}\label{eq2:lma:Bounds for Psi}
\frac{\frac{d\left(\frac{f(e^{qx})}{f^q(e^x)}\right) }{dx}}{\frac{f(e^{qx})}{f^q(e^x)}}
&=& q\left(\frac{d(\log f(e^{qx}))}{dx} - \frac{d(\log f(e^x))}{dx}\right).
\end{eqnarray}
From Eq.~\eqref{eq2:lma:Bounds for Psi}, we have $\frac{f(x^q)}{f^q(x)}$ is an decreasing function for $0 $$<$$ q $$\leq$$ 1$, and we have $\frac{f(x^q)}{f^q(x)}$ is an increasing function for $q $$\geq$$ 1$. Therefore, we have 
\begin{eqnarray}\label{eq3:lma:Bounds for Psi}
\mathcal{I} &\preceq& f^{-q_0}\left(\mathcal{Z}_n\right)f\left(\mathcal{Z}^{q_0}_n\right) 
\preceq \psi(q_0,f,\mathcal{Z}_n),
\end{eqnarray}
and, for $k=1,2,\cdots,n$, we also have
\begin{eqnarray}\label{eq4:lma:Bounds for Psi}
\mathcal{I} &\preceq& f^{-2}\left(\mathcal{Z}_k\right)f\left(\mathcal{Z}^{2}_k\right) 
\preceq \psi(2,f,\mathcal{Z}_k).
\end{eqnarray}
This Lemma is proved from Eq.~\eqref{eq3:lma:Bounds for Psi} and Eq.~\eqref{eq4:lma:Bounds for Psi} with Eq.~\eqref{eq:Bounds for Psi upper lower}.
$\hfill \Box$

The results presented in Lemma~\ref{lma:Bounds for Psi} can have us determine the tail bounds evaluation more easily in Theorem~\ref{thm:3.2} and Corollary~\ref{cor:3.3} by sacrificing precision.

\section{Tail Bounds for Bivariate Random Tensor Means Based on Majorization Ordering}\label{sec:Tail Bounds for Bivariate Random Tensors Mean Based on Majorization Ordering} 

In this section, we will derive tail bounds relations for the summation and product of eigenvalues based on majorization ordering among bivariate random tensor means.

Let $\mathbf{x} = [x_1, \dots,x_n]^{\mathrm{T}} \in \mathbb{R}^{n}, \mathbf{y} = [y_1, \dots,y_n]^{\mathrm{T}} \in \mathbb{R}^{n}$ be two vectors with the following orders among entries $x_1 \geq \cdots \geq x_n$ and $y_1 \geq \cdots \geq y_n$, \emph{weak majorization} between vectors $\mathbf{x}, \mathbf{y}$, represented by $\mathbf{x} \triangleleft_{w} \mathbf{y}$, requires the following relation for vectors $\mathbf{x}, \mathbf{y}$:
\begin{eqnarray}\label{eq:weak majorization def}
\sum\limits_{i=1}^k x_i \leq \sum\limits_{i=1}^k y_i,
\end{eqnarray}
where $k \in \{1,2,\dots,n\}$. \emph{Majorization} between vectors $\mathbf{x}, \mathbf{y}$, indicated by $\mathbf{x} \triangleleft \mathbf{y}$, needs the following relation for vectors $\mathbf{x}, \mathbf{y}$:
\begin{eqnarray}\label{eq:majorization def}
\sum\limits_{i=1}^k x_i &\leq& \sum\limits_{i=1}^k y_i,~~\mbox{for $1 \leq k < n$;} \nonumber \\
\sum\limits_{i=1}^r x_i &=& \sum\limits_{i=1}^r y_i,~~\mbox{for $k = n$.}
\end{eqnarray}

For $\mathbf{x}, \mathbf{y} \in \mathbb{R}^n_{\geq 0}$ such that  $x_1 \geq \cdots \geq x_n$ and $y_1 \geq \cdots \geq y_n$,  \emph{weak log majorization} between vectors $\mathbf{x}, \mathbf{y}$, represented by $\mathbf{x} \triangleleft_{w \log} \mathbf{y}$, needs the following relation for vectors $\mathbf{x}, \mathbf{y}$:
\begin{eqnarray}\label{eq:weak log majorization def}
\prod\limits_{i=1}^k x_i \leq \prod\limits_{i=1}^k y_i,
\end{eqnarray}
where $k \in \{1,2,\dots,n\}$, and \emph{log majorization} between vectors $\mathbf{x}, \mathbf{y}$, represented by $\mathbf{x} \triangleleft_{\log} \mathbf{y}$, requires equality for $k=n$ in Eq.~\eqref{eq:weak log majorization def}.

We need the following lemma to identify the relationships between tail bounds of different bivariate random tensor means with L\"owner ordering.
\begin{lemma}\label{lma:majorization inequality from Lowner ordering}
Given the following three random PD tensors $\mathcal{X}, \mathcal{Y}, \mathcal{Z} $$\in$$ \mathbb{C}^{I_1 \times \cdots \times I_N \times I_1 \times \cdots \times I_N}$ with the relation $\mathcal{X} $$\preceq$$ \mathcal{Y} $$\preceq$$ \mathcal{Z}$ almost surely, and their eigenvalues are arranged as 
\begin{eqnarray}\label{eq1:lma:majorization inequality from Lowner ordering}
\lambda_{1}\left(\mathcal{X}\right) \geq \lambda_{2}\left(\mathcal{X}\right) \geq \cdots \geq \lambda_{\tiny \mbox{$\prod\limits_{i=1}^N I_i$}}\left(\mathcal{X}\right), \nonumber \\
\lambda_{1}\left(\mathcal{Y}\right) \geq \lambda_{2}\left(\mathcal{Y}\right) \geq \cdots \geq \lambda_{\tiny \mbox{$\prod\limits_{i=1}^N I_i$}}\left(\mathcal{Y}\right), \nonumber \\
\lambda_{1}\left(\mathcal{Z}\right) \geq \lambda_{2}\left(\mathcal{Z}\right) \geq \cdots \geq \lambda_{\tiny \mbox{$\prod\limits_{i=1}^N I_i$}}\left(\mathcal{Z}\right).
\end{eqnarray}
Then, given any positive number $\kappa$ and $1 \leq k \leq \prod\limits_{i=1}^N I_i$, we have
\begin{eqnarray}\label{eq2:lma:majorization inequality from Lowner ordering}
\mathrm{Pr}\left(\sum\limits_{i=1}^k\lambda_{i}\left(\mathcal{X}\right) \geq \kappa\right) \leq \mathrm{Pr}\left(\sum\limits_{i=1}^k\lambda_{i}\left(\mathcal{Y}\right) \geq \kappa\right) \leq \mathrm{Pr}\left(\sum\limits_{i=1}^k\lambda_{i}\left(\mathcal{Z}\right) \geq \kappa\right),
\end{eqnarray}
and
\begin{eqnarray}\label{eq3:lma:majorization inequality from Lowner ordering}
\mathrm{Pr}\left(\prod\limits_{i=1}^k\lambda_{i}\left(\mathcal{X}\right) \geq \kappa\right) \leq \mathrm{Pr}\left(\prod\limits_{i=1}^k\lambda_{i}\left(\mathcal{Y}\right) \geq \kappa\right) \leq \mathrm{Pr}\left(\prod\limits_{i=1}^k\lambda_{i}\left(\mathcal{Z}\right) \geq \kappa\right).
\end{eqnarray}
\end{lemma}
\textbf{Proof:}
Because $\mathcal{X}$ is a Hermitian tensor, from Courant-Fisher theorem, we have 
\begin{eqnarray}
\lambda_i(\mathcal{X}) &=& \min\limits_{\bm{U}}\max\limits_{\bm{x}} \bm{x}\bm{U}\bm{X}\bm{U}^{\mathrm{H}}
\bm{x}^{\mathrm{H}},
\end{eqnarray} 
where the matrix $\bm{X}$ is the unfolded matrix from the tensor $\mathcal{X}$ according to Section 2.2 in~\cite{liang2019further}, and the matrix $\bm{U}$ runs over all $r \times \prod\limits_{j=1}^N I_j$ complex matrices satisfying $\bm{U}\bm{U}^{\mathrm{H}}$$=$$\bm{I}_{r}$. Note that $1 \leq r \leq \prod\limits_{j=1}^N I_j$. Then, for $1 $$\leq$$ i $$\leq$$ \prod\limits_{j=1}^N I_j$, we obtain the following relation:
\begin{eqnarray}\label{eq4:lma:majorization inequality from Lowner ordering}
\lambda_i(\mathcal{X}) \leq \lambda_i(\mathcal{Y}) \leq \lambda_i(\mathcal{Z}),
\end{eqnarray}
due to  $\mathcal{X} $$\preceq$$ \mathcal{Y} $$\preceq$$ \mathcal{Z}$. Since all $\lambda_i(\mathcal{X}), \lambda_i(\mathcal{Y})$ and $\lambda_i(\mathcal{Z})$ are positive numbers, this Lemma is proved from Eq.~\eqref{eq4:lma:majorization inequality from Lowner ordering}.
$\hfill \Box$

Following corollary is obtained according to Theorem~\ref{thm:3.2} to identify the relationships between tail bounds of bivariate random tensor means $\mathcal{X}^q \#_f \mathcal{Y}^q$ and $\mathcal{X} \#_f \mathcal{Y}$ for $q >0$. 

\begin{corollary}\label{thm:3.2 major}
Given two random PD tensors $\mathcal{X} $$\in$$ \mathbb{C}^{I_1 \times \cdots \times I_N \times I_1 \times \cdots \times I_N}$, $\mathcal{Y} $$\in$$ \mathbb{C}^{I_1 \times \cdots \times I_N \times I_1 \times \cdots \times I_N}$ and a PD determinstic tensor $\mathcal{C}$, if $q=2^n q_0 \geq 1$ with $1 \leq q_0 \leq 2$, we set $\mathcal{Z}_{k-1} \define \mathcal{X}^{-2^{k-2}} \mathcal{Y}^{2^{k-1}}\mathcal{X}^{-2^{k-2}}$ for $k=1,2,\cdots,n$. We assume that $\mathcal{X} $$\#_f$$ \mathcal{Y} $$\succeq$$ \mathcal{I}$ almost surely with $f $$\in$$\mbox{TMI}^{1}$. 
Then, we have 
\begin{eqnarray}\label{eq1-1:thm:3.2 major}
\mathrm{Pr}\left(\sum\limits_{i=1}^{k}\lambda_i\left(\Psi_{lower}\left(q,f,\mathcal{X},\mathcal{Y}\right)\lambda^{q-1}_{\max}\left(\mathcal{X}\#_f\mathcal{Y}\right)\mathcal{X}\#_f\mathcal{Y}\right) \geq \kappa \right)
\leq
\mathrm{Pr}\left(\sum\limits_{i=1}^{k}\lambda_i\left(\mathcal{X}^q \#_f \mathcal{Y}^q \right) \geq \kappa \right) \nonumber \\
\leq \mathrm{Pr}\left(\sum\limits_{i=1}^{k}\lambda_i\left(\Psi_{upper}\left(q,f,\mathcal{X},\mathcal{Y}\right)\lambda^{q-1}_{\min}\left(\mathcal{X}\#_f\mathcal{Y}\right)\mathcal{X}\#_f\mathcal{Y}\right) \geq \kappa \right),~~~~~~~~~~~~~~~~~~~~~~~~~~~~~~~~~~~~~~~
\end{eqnarray}
and
\begin{eqnarray}\label{eq1-2:thm:3.2 major}
\mathrm{Pr}\left(\prod\limits_{i=1}^{k}\lambda_i\left(\Psi_{lower}\left(q,f,\mathcal{X},\mathcal{Y}\right)\lambda^{q-1}_{\max}\left(\mathcal{X}\#_f\mathcal{Y}\right)\mathcal{X}\#_f\mathcal{Y}\right) \geq \kappa \right)
\leq
\mathrm{Pr}\left(\prod\limits_{i=1}^{k}\lambda_i\left(\mathcal{X}^q \#_f \mathcal{Y}^q \right) \geq \kappa \right) \nonumber \\
\leq \mathrm{Pr}\left(\prod\limits_{i=1}^{k}\lambda_i\left(\Psi_{upper}\left(q,f,\mathcal{X},\mathcal{Y}\right)\lambda^{q-1}_{\min}\left(\mathcal{X}\#_f\mathcal{Y}\right)\mathcal{X}\#_f\mathcal{Y}\right) \geq \kappa \right).~~~~~~~~~~~~~~~~~~~~~~~~~~~~~~~~~~~~~~~
\end{eqnarray}

For $0 < q  \leq1$, we have 
\begin{eqnarray}\label{eq2-1:thm:3.2 major}
\mathrm{Pr}\left(\sum\limits_{i=1}^k\lambda_i\left(\lambda_{\max}\left(f^{-q}\left(\mathcal{Z}_0\right)f\left(\mathcal{Z}_0^q\right)\right)\lambda^{q-1}_{\max}\left(\mathcal{X}\#_f\mathcal{Y}\right)\mathcal{X}\#_f\mathcal{Y}\right) \geq \kappa \right) \leq \mathrm{Pr}\left(\sum\limits_{i=1}^k\lambda_i\left(\mathcal{X}^q \#_f \mathcal{Y}^q \right) \geq \kappa \right) \nonumber \\
\leq  \mathrm{Pr}\left(\sum\limits_{i=1}^k\lambda_i\left(\lambda_{\min}\left(f^{-q}\left(\mathcal{Z}_0\right)(f\left(\mathcal{Z}_0^q\right)\right)\lambda^{q-1}_{\min}\left(\mathcal{X}\#_f\mathcal{Y}\right)\mathcal{X}\#_f\mathcal{Y} \right) \geq \kappa \right),~~~~~~~~~~~~~~~~~~~~~~~~~~~~~~~~~~~~~
\end{eqnarray}
and
\begin{eqnarray}\label{eq2-2:thm:3.2 major}
\mathrm{Pr}\left(\prod\limits_{i=1}^k\lambda_i\left(\lambda_{\max}\left(f^{-q}\left(\mathcal{Z}_0\right)f\left(\mathcal{Z}_0^q\right)\right)\lambda^{q-1}_{\max}\left(\mathcal{X}\#_f\mathcal{Y}\right)\mathcal{X}\#_f\mathcal{Y}\right) \geq \kappa \right) \leq \mathrm{Pr}\left(\prod\limits_{i=1}^k\lambda_i\left(\mathcal{X}^q \#_f \mathcal{Y}^q \right) \geq \kappa \right) \nonumber \\
\leq  \mathrm{Pr}\left(\prod\limits_{i=1}^k\lambda_i\left(\lambda_{\min}\left(f^{-q}\left(\mathcal{Z}_0\right)(f\left(\mathcal{Z}_0^q\right)\right)\lambda^{q-1}_{\min}\left(\mathcal{X}\#_f\mathcal{Y}\right)\mathcal{X}\#_f\mathcal{Y} \right) \geq \kappa \right).~~~~~~~~~~~~~~~~~~~~~~~~~~~~~~~~~~~~~
\end{eqnarray}
\end{corollary}
\textbf{Proof:}
From Eq.~\eqref{eq5-3:thm:3.2}, and the condition $q \geq 1$, we have 
\begin{eqnarray}\label{eq5-3:thm:3.2 major}
\Psi_{lower}\left(q,f,\mathcal{X},\mathcal{Y}\right)\lambda^{q-1}_{\max}\left(\mathcal{X}\#_f\mathcal{Y}\right)\mathcal{X}\#_f\mathcal{Y} &\preceq& \mathcal{X}^q \#_f \mathcal{Y}^q  \nonumber \\
&\preceq &\Psi_{upper}\left(q,f,\mathcal{X},\mathcal{Y}\right)\lambda^{q-1}_{\min}\left(\mathcal{X}\#_f\mathcal{Y}\right)\mathcal{X}\#_f\mathcal{Y}.
\end{eqnarray}
By applying Lemma~\ref{lma:majorization inequality from Lowner ordering} to Eq.~\eqref{eq5-3:thm:3.2 major}, we have Eq.~\eqref{eq1-1:thm:3.2 major} and Eq.~\eqref{eq1-2:thm:3.2 major}.

From Eq.~\eqref{eq6-4:thm:3.2}, and the condition $0 < q  \leq 1$, we have  
\begin{eqnarray}\label{eq6-4:thm:3.2 major}
\lambda_{\max}\left(f^{-q}\left(\mathcal{Z}_0\right)f\left(\mathcal{Z}_0^q\right)\right)\lambda^{q-1}_{\max}\left(\mathcal{X}\#_f\mathcal{Y}\right)\mathcal{X}\#_f\mathcal{Y} &\preceq& \mathcal{X}^q \#_f \mathcal{Y}^q  \nonumber \\
&\preceq &  \lambda_{\min}\left(f^{-q}\left(\mathcal{Z}_0\right)(f\left(\mathcal{Z}_0^q\right)\right)\lambda^{q-1}_{\min}\left(\mathcal{X}\#_f\mathcal{Y}\right)\mathcal{X}\#_f\mathcal{Y}. \nonumber \\
\end{eqnarray}
By applying Lemma~\ref{lma:majorization inequality from Lowner ordering} to Eq.~\eqref{eq6-4:thm:3.2 major}, we have Eq.~\eqref{eq2-1:thm:3.2 major} and Eq.~\eqref{eq2-2:thm:3.2 major}.
$\hfill \Box$

Following Corollary is obtained according to Theorem~\ref{thm:Prop. 3.9} to identify the relationships between tail bounds of bivariate random tensor means $\mathcal{X}^q \#_h \mathcal{Y}^q$ and $\mathcal{X} \#_h \mathcal{Y}$ for $q >0$. 

\begin{corollary}\label{Prop. 3.9 major}
Given two random PD tensors $\mathcal{X} $$\in$$ \mathbb{C}^{I_1 \times \cdots \times I_N \times I_1 \times \cdots \times I_N}$, $\mathcal{Y} $$\in$$ \mathbb{C}^{I_1 \times \cdots \times I_N \times I_1 \times \cdots \times I_N}$ and a PD determinstic tensor $\mathcal{C}$, if $q=2^n q_0 \geq 1$ with $1 \leq q_0 \leq 2$, we set $\mathcal{Z}_{k-1} \define \mathcal{X}^{-2^{k-2}} \mathcal{Y}^{2^{k-1}}\mathcal{X}^{-2^{k-2}}$ for $k=1,2,\cdots,n$. We assume that $\mathcal{X} $$\#_h$$ \mathcal{Y} $$\preceq$$ \mathcal{I}$ almost surely with $h $$\in$$\mbox{TMD}^{1}$. 
Then, we have 
\begin{eqnarray}\label{eq1-1:Prop. 3.9 major}
\mathrm{Pr}\left(\sum\limits_{i=1}^k\lambda_i\left(\Phi_{lower}\left(q,h,\mathcal{X},\mathcal{Y}\right)\lambda^{q-1}_{\max}\left(\mathcal{X}\#_h\mathcal{Y}\right)\mathcal{X}\#_h\mathcal{Y}\right) \geq \kappa \right) \leq \mathrm{Pr}\left(\sum\limits_{i=1}^k\lambda_i\left(\mathcal{X}^q \#_h \mathcal{Y}^q \right) \geq \kappa \right)  \nonumber \\
\leq \mathrm{Pr}\left(\sum\limits_{i=1}^k\lambda_i\left(\Phi_{upper}\left(q,h,\mathcal{X},\mathcal{Y}\right)\lambda^{q-1}_{\min}\left(\mathcal{X}\#_h\mathcal{Y}\right)\mathcal{X}\#_h\mathcal{Y}\right) \geq \kappa \right),~~~~~~~~~~~~~~~~~~~~~~~~~~~~~~~~~~~~~~
\end{eqnarray}
and
\begin{eqnarray}\label{eq1-2:Prop. 3.9 major}
\mathrm{Pr}\left(\prod\limits_{i=1}^k\lambda_i\left(\Phi_{lower}\left(q,h,\mathcal{X},\mathcal{Y}\right)\lambda^{q-1}_{\max}\left(\mathcal{X}\#_h\mathcal{Y}\right)\mathcal{X}\#_h\mathcal{Y}\right) \geq \kappa \right) \leq \mathrm{Pr}\left(\prod\limits_{i=1}^k\lambda_i\left(\mathcal{X}^q \#_h \mathcal{Y}^q \right) \geq \kappa \right)  \nonumber \\
\leq \mathrm{Pr}\left(\prod\limits_{i=1}^k\lambda_i\left(\Phi_{upper}\left(q,h,\mathcal{X},\mathcal{Y}\right)\lambda^{q-1}_{\min}\left(\mathcal{X}\#_h\mathcal{Y}\right)\mathcal{X}\#_h\mathcal{Y}\right) \geq \kappa \right).~~~~~~~~~~~~~~~~~~~~~~~~~~~~~~~~~~~~~
\end{eqnarray}

For $0 < q  \leq1$, we have 
\begin{eqnarray}\label{eq2-1:Prop. 3.9 major}
\mathrm{Pr}\left(\sum\limits_{i=1}^k\lambda_i\left(\lambda_{\min}\left(h^{-q}\left(\mathcal{Z}_0\right)h\left(\mathcal{Z}_0^q\right)\right)\lambda^{q-1}_{\max}\left(\mathcal{X}\#_h\mathcal{Y}\right)\mathcal{X}\#_h\mathcal{Y}\right) \geq \kappa \right) \leq \mathrm{Pr}\left(\sum\limits_{i=1}^k\lambda_i\left(\mathcal{X}^q \#_h \mathcal{Y}^q \right) \geq \kappa \right) \nonumber \\
\leq  \mathrm{Pr}\left(\sum\limits_{i=1}^k\lambda_i\left(\lambda_{\max}\left(h^{-q}\left(\mathcal{Z}_0\right)(h\left(\mathcal{Z}_0^q\right)\right)\lambda^{q-1}_{\min}\left(\mathcal{X}\#_h\mathcal{Y}\right)\mathcal{X}\#_h\mathcal{Y}\right) \geq \kappa \right),~~~~~~~~~~~~~~~~~~~~~~~~~~~~~~~~~~~~~
\end{eqnarray}
and
\begin{eqnarray}\label{eq2-2:Prop. 3.9 major}
\mathrm{Pr}\left(\prod\limits_{i=1}^k\lambda_i\left(\lambda_{\min}\left(h^{-q}\left(\mathcal{Z}_0\right)h\left(\mathcal{Z}_0^q\right)\right)\lambda^{q-1}_{\max}\left(\mathcal{X}\#_h\mathcal{Y}\right)\mathcal{X}\#_h\mathcal{Y}\right) \geq \kappa \right) \leq \mathrm{Pr}\left(\prod\limits_{i=1}^k\lambda_i\left(\mathcal{X}^q \#_h \mathcal{Y}^q \right) \geq \kappa \right) \nonumber \\
\leq  \mathrm{Pr}\left(\prod\limits_{i=1}^k\lambda_i\left(\lambda_{\max}\left(h^{-q}\left(\mathcal{Z}_0\right)(h\left(\mathcal{Z}_0^q\right)\right)\lambda^{q-1}_{\min}\left(\mathcal{X}\#_h\mathcal{Y}\right)\mathcal{X}\#_h\mathcal{Y}\right) \geq \kappa \right).~~~~~~~~~~~~~~~~~~~~~~~~~~~~~~~~~~~~~
\end{eqnarray}
\end{corollary}
\textbf{Proof:}
From Eq.~\eqref{eq5-3:Prop. 3.9}, and the condition $q \geq 1$, we have 
\begin{eqnarray}\label{eq5-3:Prop. 3.9 major}
\Phi_{lower}\left(q,h,\mathcal{X},\mathcal{Y}\right)\lambda^{q-1}_{\max}\left(\mathcal{X}\#_h\mathcal{Y}\right)\mathcal{X}\#_h\mathcal{Y} &\preceq& \mathcal{X}^q \#_h \mathcal{Y}^q  \nonumber \\
&\preceq &\Phi_{upper}\left(q,h,\mathcal{X},\mathcal{Y}\right)\lambda^{q-1}_{\min}\left(\mathcal{X}\#_h\mathcal{Y}\right)\mathcal{X}\#_h\mathcal{Y}.
\end{eqnarray}
By applying Lemma~\ref{lma:majorization inequality from Lowner ordering} to Eq.~\eqref{eq5-3:Prop. 3.9 major}, we have Eq.~\eqref{eq1-1:Prop. 3.9 major} and Eq.~\eqref{eq1-2:Prop. 3.9 major}.

From Eq.~\eqref{eq6-4:Prop. 3.9}, and the condition $0 < q  \leq 1$, we have 
\begin{eqnarray}\label{eq6-4:Prop. 3.9 major}
\lambda_{\min}\left(h^{-q}\left(\mathcal{Z}_0\right)h\left(\mathcal{Z}_0^q\right)\right)\lambda^{q-1}_{\max}\left(\mathcal{X}\#_h\mathcal{Y}\right)\mathcal{X}\#_h\mathcal{Y} &\preceq& \mathcal{X}^q \#_h \mathcal{Y}^q  \nonumber \\
&\preceq &  \lambda_{\max}\left(h^{-q}\left(\mathcal{Z}_0\right)(h\left(\mathcal{Z}_0^q\right)\right)\lambda^{q-1}_{\min}\left(\mathcal{X}\#_h\mathcal{Y}\right)\mathcal{X}\#_h\mathcal{Y}. \nonumber \\
\end{eqnarray}
By applying Lemma~\ref{lma:majorization inequality from Lowner ordering} to Eq.~\eqref{eq6-4:Prop. 3.9 major}, we have Eq.~\eqref{eq2-1:Prop. 3.9 major} and Eq.~\eqref{eq2-2:Prop. 3.9 major}.
$\hfill \Box$

Following corollary is obtained according to Theorem~\ref{thm:Prop. 3.10} to identify the relationships between tail bounds of bivariate random tensor means $\mathcal{X}^q \#_f \mathcal{Y}^q$ and $\mathcal{X} \#_f \mathcal{Y}$ for $q >0$. 

\begin{corollary}\label{thm:Prop. 3.10 major}
Given two random PD tensors $\mathcal{X} $$\in$$ \mathbb{C}^{I_1 \times \cdots \times I_N \times I_1 \times \cdots \times I_N}$, $\mathcal{Y} $$\in$$ \mathbb{C}^{I_1 \times \cdots \times I_N \times I_1 \times \cdots \times I_N}$, and a PD determinstic tensor $\mathcal{C} $$\in$$ \mathbb{C}^{I_1 \times \cdots \times I_N \times I_1 \times \cdots \times I_N}$, we will set $\mathcal{Z} \define \mathcal{X}^{1/2}\star_N\mathcal{Y}^{-1}\star_N\mathcal{X}^{1/2}$. Let $g $$\in$$\mbox{TC}^{1}$, if $\mathcal{X} \#_g \mathcal{Y} \preceq \mathcal{I}$ almost surely and $q \geq 1$. Then, we have 
\begin{eqnarray}\label{eq1-1:thm:Prop. 3.10 major}
\mathrm{Pr}\left(\sum\limits_{i=1}^k\lambda_i\left(\mathcal{X}^q \#_g \mathcal{Y}^q  \right) \geq \kappa \right) &\leq& \mathrm{Pr}\Bigg(\sum\limits_{i=1}^k\lambda_i\Big(\mathrm{K}\left(\lambda^{-1}_{\max}\left(\mathcal{X}\right),\lambda^{-1}_{\min}\left(\mathcal{X}\right),q-1\right)\lambda^{1-q}_{\min}\left(\mathcal{X}\#_{g}\mathcal{Y}\right) \nonumber \\
&  & \lambda_{\max}\left(g^{-q}(\mathcal{Z})g(\mathcal{Z}^q)\right)\mathrm{K}\left(\lambda^{-1}_{\max}\left(\mathcal{X}\right),\lambda^{-1}_{\min}\left(\mathcal{X}\right),2q-1\right)\mathcal{I}\Big) \geq \kappa \Bigg),\nonumber \\
\end{eqnarray}
and
\begin{eqnarray}\label{eq1-2:thm:Prop. 3.10 major}
\mathrm{Pr}\left(\prod\limits_{i=1}^k\lambda_i\left(\mathcal{X}^q \#_g \mathcal{Y}^q  \right) \geq \kappa \right) &\leq& \mathrm{Pr}\Bigg(\prod\limits_{i=1}^k\lambda_i\Big(\mathrm{K}\left(\lambda^{-1}_{\max}\left(\mathcal{X}\right),\lambda^{-1}_{\min}\left(\mathcal{X}\right),q-1\right)\lambda^{1-q}_{\min}\left(\mathcal{X}\#_{g}\mathcal{Y}\right) \nonumber \\
&  & \lambda_{\max}\left(g^{-q}(\mathcal{Z})g(\mathcal{Z}^q)\right)\mathrm{K}\left(\lambda^{-1}_{\max}\left(\mathcal{X}\right),\lambda^{-1}_{\min}\left(\mathcal{X}\right),2q-1\right)\mathcal{I}\Big) \geq \kappa \Bigg)\nonumber \\
\end{eqnarray}

Moreover, if $\mathcal{X} \#_g \mathcal{Y} \succeq \mathcal{I}$ almost surely, we have 
\begin{eqnarray}\label{eq2-1:thm:Prop. 3.10 major}
\mathrm{Pr}\left(\sum\limits_{i=1}^k\lambda_i\left(\mathcal{X}^q \#_g \mathcal{Y}^q  \right) \geq \kappa \right) &\geq& \mathrm{Pr}\Bigg(\sum\limits_{i=1}^k\lambda_i\Big(\lambda^{1-q}_{\min}\left(\mathcal{X}\#_{g}\mathcal{Y}\right)\lambda_{\max}\left(g^{-q}(\mathcal{Z})g(\mathcal{Z}^q)\right) \nonumber \\
&  & \mathrm{K}^{-1}\left(\lambda^{-1}_{\max}\left(\mathcal{X}\right),\lambda^{-1}_{\min}\left(\mathcal{X}\right),2q-1\right)\mathcal{I}\Big) \geq \kappa \Bigg),
\end{eqnarray}
and
\begin{eqnarray}\label{eq2-2:thm:Prop. 3.10 major}
\mathrm{Pr}\left(\prod\limits_{i=1}^k\lambda_i\left(\mathcal{X}^q \#_g \mathcal{Y}^q  \right) \geq \kappa \right) &\geq& \mathrm{Pr}\Bigg(\prod\limits_{i=1}^k\lambda_i\Big(\lambda^{1-q}_{\min}\left(\mathcal{X}\#_{g}\mathcal{Y}\right)\lambda_{\max}\left(g^{-q}(\mathcal{Z})g(\mathcal{Z}^q)\right) \nonumber \\
&  & \mathrm{K}^{-1}\left(\lambda^{-1}_{\max}\left(\mathcal{X}\right),\lambda^{-1}_{\min}\left(\mathcal{X}\right),2q-1\right)\mathcal{I}\Big) \geq \kappa \Bigg).
\end{eqnarray}
\end{corollary}
\textbf{Proof:}
From Eq.~\eqref{eq4-1:Prop. 3.10}, and the condition $q \geq 1$, we have 
\begin{eqnarray}\label{eq4-1:Prop. 3.10 major}
\mathcal{X}^q \#_g \mathcal{Y}^q &\preceq& \mathrm{K}\left(\lambda^{-1}_{\max}\left(\mathcal{X}\right),\lambda^{-1}_{\min}\left(\mathcal{X}\right),q-1\right)\lambda^{1-q}_{\min}\left(\mathcal{X}\#_{g}\mathcal{Y}\right)
\lambda_{\max}\left(g^{-q}(\mathcal{Z})g(\mathcal{Z}^q)\right) \nonumber \\
&  &\mathrm{K}\left(\lambda^{-1}_{\max}\left(\mathcal{X}\right),\lambda^{-1}_{\min}\left(\mathcal{X}\right),2q-1\right)\mathcal{I}.
\end{eqnarray}
By applying Lemma~\ref{lma:majorization inequality from Lowner ordering} to Eq.~\eqref{eq4-1:Prop. 3.10 major}, we have Eq.~\eqref{eq1-1:thm:Prop. 3.10 major} and Eq.~\eqref{eq1-2:thm:Prop. 3.10 major}.

From Eq.~\eqref{eq7-1:Prop. 3.10}, and the condition $0 < q  \leq 1$, we have 
\begin{eqnarray}\label{eq7-1:Prop. 3.10 major}
\mathcal{X}^q \#_g \mathcal{Y}^q &\succeq& \lambda^{1-q}_{\min}\left(\mathcal{X}\#_{g}\mathcal{Y}\right)
\lambda_{\max}\left(g^{-q}(\mathcal{Z})g(\mathcal{Z}^q)\right) \mathrm{K}^{-1}\left(\lambda^{-1}_{\max}\left(\mathcal{X}\right),\lambda^{-1}_{\min}\left(\mathcal{X}\right),2q-1\right)\mathcal{I},
\end{eqnarray}
By applying Lemma~\ref{lma:majorization inequality from Lowner ordering} to Eq.~\eqref{eq7-1:Prop. 3.10 major}, we have Eq.~\eqref{eq2-1:thm:Prop. 3.10 major} and Eq.~\eqref{eq2-2:thm:Prop. 3.10 major}.
$\hfill \Box$

\bibliographystyle{IEEETran}
\bibliography{AHType_2Argu_Bib}

\begin{thebibliography}{10}
\providecommand{\url}[1]{#1}
\csname url@samestyle\endcsname
\providecommand{\newblock}{\relax}
\providecommand{\bibinfo}[2]{#2}
\providecommand{\BIBentrySTDinterwordspacing}{\spaceskip=0pt\relax}
\providecommand{\BIBentryALTinterwordstretchfactor}{4}
\providecommand{\BIBentryALTinterwordspacing}{\spaceskip=\fontdimen2\font plus
\BIBentryALTinterwordstretchfactor\fontdimen3\font minus
  \fontdimen4\font\relax}
\providecommand{\BIBforeignlanguage}[2]{{%
\expandafter\ifx\csname l@#1\endcsname\relax
\typeout{** WARNING: IEEEtran.bst: No hyphenation pattern has been}%
\typeout{** loaded for the language `#1'. Using the pattern for}%
\typeout{** the default language instead.}%
\else
\language=\csname l@#1\endcsname
\fi
#2}}
\providecommand{\BIBdecl}{\relax}
\BIBdecl

\bibitem{ouerfelli2022random}
M.~Ouerfelli, M.~Tamaazousti, and V.~Rivasseau, ``Random tensor theory for
  tensor decomposition,'' in \emph{Proceedings of the AAAI Conference on
  Artificial Intelligence}, vol.~36, no.~7, 2022, pp. 7913--7921.

\bibitem{evnin2021melonic}
O.~Evnin, ``Melonic dominance and the largest eigenvalue of a large random
  tensor,'' \emph{Letters in Mathematical Physics}, vol. 111, no.~3, p.~66,
  2021.

\bibitem{tropp2015introduction}
J.~A. Tropp \emph{et~al.}, ``An introduction to matrix concentration
  inequalities,'' \emph{Foundations and Trends in Machine Learning}, vol.~8,
  no. 1-2, pp. 1--230, 2015.

\bibitem{chang2022TWF}
S.~Y. Chang and H.-C. Wu, ``Tensor wiener filter,'' \emph{IEEE Transactions on
  Signal Processing}, vol.~70, pp. 410--422, 2022.

\bibitem{chang2022randouble}
S.~Y. Chang, ``Random double tensors integrals,'' \emph{arXiv preprint
  arXiv:2204.01927}, 2022.

\bibitem{chang2022tPII}
S.~Y. Chang and Y.~Wei, ``T-product tensors—part ii: tail bounds for sums of
  random t-product tensors,'' \emph{Computational and Applied Mathematics},
  vol.~41, no.~3, p.~99, 2022.

\bibitem{chang2022tPI}
------, ``T-square tensors—part i: inequalities,'' \emph{Computational and
  Applied Mathematics}, vol.~41, no.~1, p.~62, 2022.

\bibitem{chang2022generaltail}
------, ``General tail bounds for random tensors summation: majorization
  approach,'' \emph{Journal of Computational and Applied Mathematics}, vol.
  416, p. 114533, 2022.

\bibitem{chang2022TKF}
S.~Y. Chang and H.-C. Wu, ``Tensor kalman filter and its applications,''
  \emph{IEEE Transactions on Knowledge and Data Engineering}, 2022.

\bibitem{chang2022tensorq}
------, ``Tensor quantization: High-dimensional data compression,'' \emph{IEEE
  Transactions on Circuits and Systems for Video Technology}, vol.~32, no.~8,
  pp. 5566--5580, 2022.

\bibitem{kubo1980means}
F.~Kubo and T.~Ando, ``Means of positive linear operators,''
  \emph{Mathematische Annalen}, vol. 246, pp. 205--224, 1980.

\bibitem{ando1994log}
T.~Ando and F.~Hiai, ``Log majorization and complementary {G}olden-{T}hompson
  type inequalities,'' \emph{Linear Algebra and its Applications}, vol. 197,
  pp. 113--131, 1994.

\bibitem{wada2018does}
S.~Wada, ``When does {A}ndo--{H}iai inequality hold?'' \emph{Linear Algebra and
  its Applications}, vol. 540, pp. 234--243, 2018.

\bibitem{hiai2020ando}
F.~Hiai, Y.~Seo, and S.~Wada, ``{A}ndo--{H}iai-type inequalities for operator
  means and operator perspectives,'' \emph{International Journal of
  Mathematics}, vol.~31, no.~01, p. 2050007, 2020.

\bibitem{effros2014non}
E.~Effros and F.~Hansen, ``Non-commutative perspectives,'' \emph{Annals of
  Functional Analysis}, vol.~5, no.~2, pp. 74--79, 2014.

\bibitem{ni2019hermitian}
G.~Ni, ``Hermitian tensor and quantum mixed state,'' \emph{arXiv preprint
  arXiv:1902.02640}, 2019.

\bibitem{furuta2005mond}
T.~Furuta and J.~M. Hot, ``{M}ond-{P}ecaric method in operator inequalities,''
  \emph{Inequalities for bounded selfadjoint operators on a Hilbert space,
  Element, Zagreb}, 2005.

\bibitem{liang2019further}
M.~Liang and B.~Zheng, ``Further results on {M}oore-{P}enrose inverses of
  tensors with application to tensor nearness problems,'' \emph{Computers and
  Mathematics with Applications}, vol.~77, no.~5, pp. 1282--1293, March 2019.

\end{thebibliography}

\end{document}